\documentclass[10pt,a4paper]{article}
\usepackage[a4paper]{geometry}
\usepackage{amssymb,latexsym,amsmath,amsfonts,amsthm}
\usepackage{graphicx}

\usepackage{epsfig}
\usepackage{comment}
\usepackage{subfigure}
\usepackage{graphicx,ebezier}
\usepackage{overpic}

\DeclareMathOperator{\diag}{diag}

\newcommand{\lam}{\lambda}

\newcommand{\crit}{\textrm{crit}}

\renewcommand{\Re}{\mathrm{Re}\,}
\renewcommand{\Im}{\mathrm{Im}\,}
\newcommand{\ds}{\displaystyle}

\newcommand{\Ai}{\mathrm{Ai}}
\newcommand{\vece}{\mathbf{e}}
\renewcommand{\vec}{\mathbf}
\newcommand{\ccal}{\mathcal{C}}
\newcommand{\supp}{\mathrm{supp}}

\newtheorem{theorem}{Theorem}[section]
\newtheorem{lemma}[theorem]{Lemma}
\newtheorem{proposition}[theorem]{Proposition}

\newtheorem{example}{Example}
\newtheorem{corollary}[theorem]{Corollary}
\newtheorem{algorithm}[theorem]{Algorithm}
\newtheorem{rhp}[theorem]{RH problem}
\newtheorem{assumption}[theorem]{Assumption}

\theoremstyle{definition}
\newtheorem{definition}[theorem]{Definition}

\theoremstyle{remark}

\newtheorem{remark}[theorem]{Remark}

\numberwithin{equation}{section}

\title{A graph-based equilibrium problem for the limiting distribution
of non-intersecting Brownian motions at low temperature}

\author{Steven Delvaux\footnotemark[1],\quad Arno B.~J.~Kuijlaars\footnotemark[1]}
\date{}

\begin{document}

\maketitle
\renewcommand{\thefootnote}{\fnsymbol{footnote}}
\footnotetext[1]{Department of Mathematics, Katholieke
Universiteit Leuven, Celestijnenlaan 200B, B-3001 Leuven, Belgium.
email:
\{steven.delvaux,arno.kuijlaars\}\symbol{'100}wis.kuleuven.be. \\
The first author is a Postdoctoral Fellow of the Fund for Scientific Research -
Flanders (Belgium). \\
The work of the second author is supported by FWO-Flanders project
G.0427.09, by K.U. Leuven research grant OT/08/33, by
the Belgian Interuniversity Attraction Pole P06/02, by the
European Science Foundation Program MISGAM, and by grant
MTM2008-06689-C02-01 of the Spanish
Ministry of Science and Innovation.}

\begin{abstract}
We consider $n$ non-intersecting Brownian motion paths with $p$ prescribed
starting positions at time $t=0$ and $q$ prescribed ending positions at time $t=1$.
The positions of the paths at any intermediate time are a determinantal
point process, which in the case $p=1$ is equivalent to the eigenvalue
distribution of a random matrix from the Gaussian unitary ensemble with external source.
For general $p$ and $q$, we show that
if a temperature parameter is sufficiently small, then the distribution
of the Brownian paths is characterized in the large $n$ limit
by a vector equilibrium problem with an interaction matrix that is based
on a bipartite planar graph.
Our proof is based on a steepest descent
analysis of an associated $(p+q) \times (p+q)$ matrix valued
Riemann-Hilbert problem whose solution is built out of multiple orthogonal
polynomials. A new feature of the steepest descent analysis is a systematic
opening of a large number of global lenses.


\end{abstract}

\section{Introduction}
\label{section:introduction}

\begin{figure}[t]
\begin{center}\vspace{-5mm}
\includegraphics[scale=0.55]{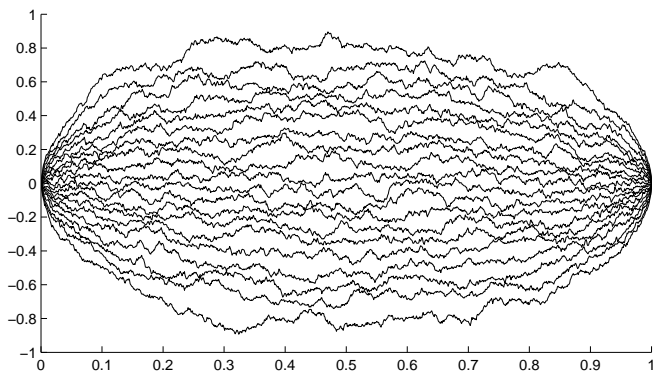}
\end{center}\vspace{-5mm}
\caption{Non-intersecting Brownian motions with one starting and one ending
point. As the number of paths tends to infinity, the paths fill out an ellipse
in the $tx$-plane.} \label{figwigner}
\end{figure}

This paper deals with non-intersecting one-dimensional Brownian motions with
prescribed starting and ending positions. This model has already been discussed
in various regimes. For the case of one starting point and one ending point it
is known that the positions of the paths at any intermediate time have the same
distribution (up to trivial scaling) as the eigenvalues of a Gaussian Unitary
Ensemble from random matrix theory \cite{Dyson}. Moreover, as the number of
paths tends to infinity and after appropriate scaling, the paths fill out an
ellipse in the $tx$-plane, see Figure \ref{figwigner}.

\begin{figure}[t]
\begin{center}
\includegraphics[scale=0.55]{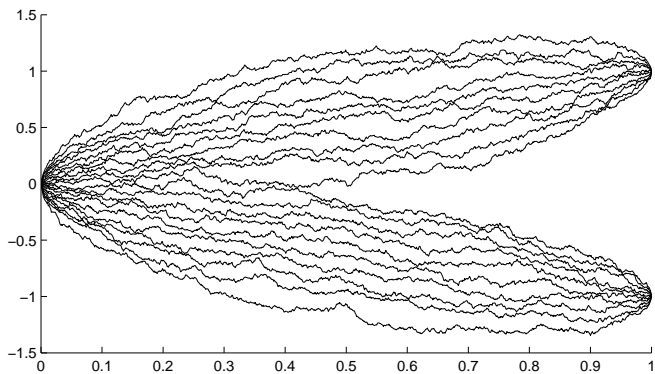}
\end{center}\vspace{-5mm}
\caption{Non-intersecting Brownian motions with one starting and two ending
positions.} \label{figwigner2}
\end{figure}

In the case of one starting point and two or more ending points the positions
of the paths have the same distribution as the eigenvalues of a Gaussian
Unitary Ensemble with external source. This model is described by multiple
Hermite polynomials. As the number of paths tends to infinity, the paths fill
out a more complicated region whose boundary has cusp points. The limiting
distributions can be computed in terms of an algebraic curve known as
\emph{Pastur's equation} \cite{ABK,BK2,BK3,KKPS,Ora,Pastur}. See Figure
\ref{figwigner2} for an illustration of the case of two ending points.

The case of one ending point and two or more starting points is
equivalent due to the time reversal symmetry in the model.

Much less is known for the general case of $p \geq 2$ starting points and $q
\geq 2$ ending points. For example, it  is not known whether there exists an
underlying random matrix model for this case. What is known is that the model
is described by multiple Hermite polynomials of mixed type \cite{DK2} which
have a characterization in terms of a $(p + q) \times (p+q)$ matrix valued
Riemann Hilbert problem. Calculations for the limiting distributions of paths
in the large $n$ limit were done for very specific cases with $p=q=2$ in
\cite{DKV,DelKui} based on the spectral curve (analogue of the Pastur's
equation) that could be computed in these special cases, see also the related
works \cite{AFvM,Mcl}.

It is the goal of this paper to study the case of general $p \geq 2$
and $q \geq 2$ with methods from potential theory, more precisely
with vector equilibrium problems with external fields.

Equilibrium problems with external fields were developed in approximation
theory in the context of orthogonal polynomials, Pad\'e approximation, and
polynomial approximation with varying weights, see \cite{SaffTotik,Totik}. They
are also a powerful tool in the study of unitary random matrix ensembles
\cite{Dei,DKM}. Vector equilibrium problems were studied in the context of
Hermite-Pad\'e approximation and the associated multiple orthogonal polynomials
\cite{Apt,NS}.

\begin{figure}[t]
\begin{center}\vspace{-5mm}
\includegraphics[scale=0.55]{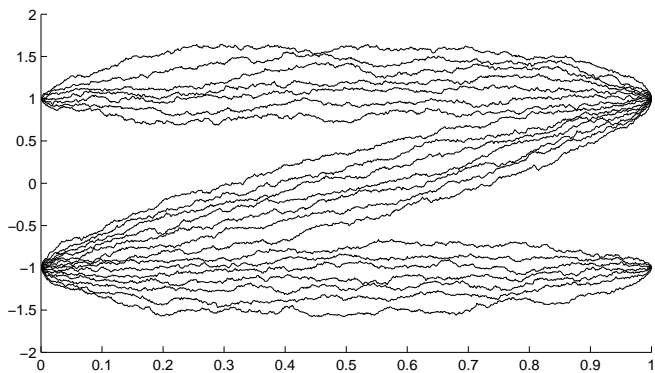}
\end{center}\vspace{-5mm}
\caption{Non-intersecting Brownian motions with two starting and two ending
positions. The starting and ending positions are sufficiently far apart so that
around time $t=1/2$, there are three groups of paths in the large $n$ limit.}
\label{figwigner3}
\end{figure}

For $p = q = 2$ we are considering a situation such as the one shown
in Figure~\ref{figwigner3}, where a certain fraction of the
paths starts in each of the two starting points, and
ends at each of the two ending points. As the number
of paths increases we see the following situation. For
small time the paths are in two separate groups that emanate
from the two starting positions. At a certain time one
of the groups splits into two, leading to a situation
of three separate groups of paths. Then at a later time
two of the groups come together and we end up with two
groups that end at the two ending points.

Our results will deal with the situation at times where there
are three groups of paths, or for general $p$ and $q$, where
there are the maximal number (namely $p+q-1$) of groups of paths.

There is an alternative possible scenario in which the two groups
of paths first merge into one group and later split again into
two groups of paths. This will happen if the starting and ending
positions are sufficiently close to each other. The first
scenario happens if the starting and ending positions are relatively
far away from each other. Below we will actually distinguish  the two
scenarios in terms of a temperature parameter $T$ so that for small $T$ we
have the situation with the three groups of paths.

\section{Statement of results} \label{subsection:statementofresults}
\subsection{Assumptions}
\label{subsection:brownintro}

Let $p \geq 2$ and $q \geq 2$. We fix $p$ starting  points $a_1, \ldots, a_p$
which we assume to be ordered as
\begin{equation} \label{apoints}
    a_1> a_2 > \cdots > a_p,
    \end{equation}
and $q$ ending points $b_1, \ldots, b_q$ with
\begin{equation} \label{bpoints}
    b_1> b_2 > \cdots > b_q.
    \end{equation}

For a given (large) $n$ we consider $n$ non-intersecting Brownian motion paths
and we assume that $n_k$ of the paths start at $a_k$ and that $m_l$ of the
paths end at $b_l$ for $k=1, \ldots, p$ and $l = 1, \ldots, q$. Thus
\[ \sum_{k=1}^p n_k = \sum_{l=1}^q m_l  = n. \]
Since the paths are non-intersecting, the numbers $n_k$ and $m_l$ also
determine for each $k=1, \ldots, p$ and $l= 1, \ldots, q$, the number $n_{k,l}$
of paths that start at $a_k$ and end at $b_l$.
We call the fractions
\begin{equation} \label{transitionnumbers}
    t_{k,l}^{(n)} = \frac{n_{k,l}}{n}
    \end{equation}
the finite $n$ \emph{transition numbers}.
Note that
\begin{equation} \label{transitiontkl}
    t_{k,l}^{(n)} \geq 0, \qquad \sum_{k=1}^p \sum_{l=1}^q t_{k,l}^{(n)} = 1.
    \end{equation}

As $n \to \infty$, we assume that the finite $n$ transition numbers
have limits
\begin{equation} \label{limittransition}
    t_{k,l} = \lim_{n \to \infty} t_{k,l}^{(n)}
    \end{equation}
which are the \emph{limiting transition numbers}.
It is convenient to arrange the (finite $n$ and limiting) transition
numbers into  $p \times q$ matrices
\[ \left( t_{k,l}^{(n)} \right)_{k=1, \ldots, p, l = 1, \ldots, q},
    \qquad \left( t_{k,l} \right)_{k=1, \ldots, p, l = 1, \ldots, q}. \]
To avoid degenerate cases, we assume that each row and column of the matrix
$\left( t_{k,l} \right)_{k=1, \ldots, p, l = 1, \ldots, q}$ has at least one
non-zero entry.

The assumption that the paths are non-intersecting puts a number of
constraints on the numbers $n_{k,l}$ and on the limiting transition
numbers $t_{k,l}$. Indeed, not all $a_k$ can be connected to all $b_l$
and certain transition numbers must be zero.
The constraints on the transition numbers are easy to visualize
in terms of a weighted bipartite graph
\begin{equation}\label{defgraph}
    G = (V,E,t),
    \end{equation}
with vertices
\begin{equation*}
    V = \{a_1, \ldots, a_p \} \sqcup \{b_1, \ldots, b_q \},
        \qquad \text{(disjoint union)},
\end{equation*}
edges
\begin{equation*}
    E = \{ (a_k,b_l) \in V \times V \mid t_{k,l} > 0 \},
\end{equation*}
and a weight function
\begin{equation}\label{deftfunction}
    t: E \to (0,1] :  \, (a_k,b_l) \mapsto t_{k,l}.
\end{equation}

\begin{example} \label{example1}
The graph $G = (V,E,t)$ associated with Figure~\ref{figwigner3} is shown in
Figure~\ref{fig:graph1}. The graph has four vertices and three edges, each of them
with weight $1/3$.
\begin{figure}[htbp]
\begin{center}\vspace{-1mm}
         \includegraphics[scale=0.6]{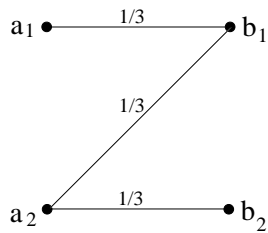}
\end{center}\vspace{-5mm}
\caption{The graph associated with Figure~\ref{figwigner3}.}
\label{fig:graph1}
\end{figure}
\end{example}

\begin{example} \label{example2}
For a more complicated example we consider a situation with $p=2$ starting
points and $q=4$ ending points as in Figure \ref{figwigner4}.

\begin{figure}[t]
\begin{center}\vspace{-1mm}
\includegraphics[scale=0.55]{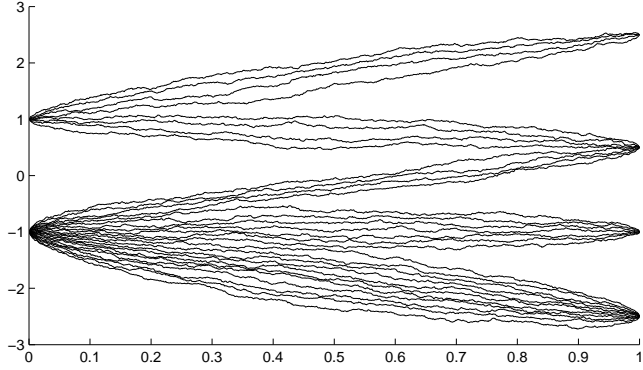}
\end{center}\vspace{-5mm}
\caption{Non-intersecting Brownian motions with two starting points and
four ending points, as discussed in Example \ref{example2}.}
\label{figwigner4}
\end{figure}

The matrix of transition numbers is
\begin{align}\label{transitionvb:matrix2}
    \begin{pmatrix}t_{k,l}\end{pmatrix}_{k=1,2,\ l=1,\ldots,4} =
    \begin{pmatrix} 4/30 & 4/30 & 0 & 0 \\
    0 & 4/30 & 7/30 & 11/30 \end{pmatrix}
\end{align}
and the graph $G$ associated with \eqref{transitionvb:matrix2} is shown in
Figure~\ref{fig:graph2}.

\begin{figure}[htbp]
\begin{center}\vspace{-1mm}
        \includegraphics[scale=0.6]{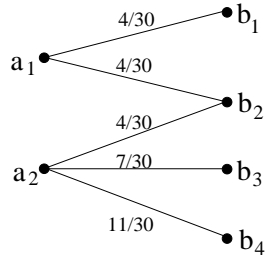}
\end{center}\vspace{-5mm}
\caption{The graph associated with the transition numbers \eqref{transitionvb:matrix2}.}
\label{fig:graph2}
\end{figure}
\end{example}

The constraints on the transition numbers are contained in
the following obvious result that we state without proof.
\begin{proposition}\label{lemma:forest}
The graph $G$ has the following properties:
\begin{enumerate}
\item[\rm (a)] $G$ has at most $p+q-1$ edges. For each $i = 1, \ldots, p+q-1$, there is at most
one non-vanishing transition number $t_{k,l}$ with $ k + l - 1 = i$.
\item[\rm (b)] $G$ is a connected graph if and only if the number of edges is equal to $p+q-1$.
\item[\rm (c)] $G$ has no cycles (and so $G$ is a tree if $G$ is connected).
\end{enumerate}
\end{proposition}

In \cite{DKV,DelKui} the special case of transition numbers
\[ \left( t_{k,l} \right)_{k=1,2, \ l = 1,2} =  \begin{pmatrix} 1/2 & 0 \\ 0 & 1/2 \end{pmatrix} \]
was considered. This example leads to a non-connected graph.

In the case of a connected graph the  structure of the matrix of
transition numbers $(t_{k,l})$ is easy to describe.

\begin{proposition} \label{lemma:latticepath}
Suppose that the graph $G$ is connected.
Then the non-zero entries of the matrix
$\begin{pmatrix}t_{k,l}\end{pmatrix}_{k,l}$ are situated on a \emph{right-down path}
starting at the top left entry $(1,1)$ and ending at the bottom right entry $(p,q)$.
The steps in the path are either by one unit to the right (a right step)
or one unit down (a down step).
\end{proposition}

In Examples \ref{example1} and \ref{example2} we have
\begin{align*}
\begin{pmatrix} \times & 0 \\
\times & \times
\end{pmatrix},\qquad
\begin{pmatrix} \times & \times & 0 & 0 \\
0 & \times & \times & \times
\end{pmatrix},
\end{align*}
respectively.

In this paper we consider only connected graphs, and we will make the
following assumption.
\begin{assumption} \label{convention:connected}
We assume that the graph $G$ is connected.  That is, we assume that $G$ has
$p+q-1$ edges, and for every $i=1, \ldots, p+q-1$ there is exactly one
non-vanishing transition number $t_{k,l}$ with $k+l-1 = i$, and we define
\[ k(i) = k, \qquad l(i) = l, \qquad \text{if } t_{k,l} > 0, \quad \text{and} \quad k+l-1 = i. \]
\end{assumption}

It follows from the assumption and from \eqref{limittransition} that
also for large enough $n$, the finite $n$ transition numbers
have the same non-zero pattern. Thus
$t_{k,l}^{(n)} > 0$ if and only if $t_{k,l} > 0$.

In what follows we follow the convention that $i$ labels the edges $E$
of the graph, and so we identify $i$ with the edge
\[ (a_{k(i)}, b_{l(i)}) \]
of the graph.

\subsection{Non-intersecting Brownian motions}

We consider Brownian motions having  transition probability density
\begin{equation} \label{transitionprob}
  P(t,x,y) = \frac{1}{\sqrt{2\pi  t \sigma^2}} e^{- \frac{1}{2t \sigma^2}(x-y)^2}
  \end{equation}
whose overall variance
\begin{equation} \label{variance}
    \sigma^2 = \frac{T}{n}
    \end{equation}
is proportional to $1/n$ where $n$ is the number of Brownian paths.
We interpret the proportionality constant $T >0$
as a temperature variable.

In this paper we consider small temperature $T$. We show
that for small $T$ the paths at time $t$ have a limiting mean
distribution that is characterized by a vector equilibrium
problem. In the first theorem we state the existence of a limiting
mean distribution.

\begin{theorem}\label{theorem:limitingdistribution}
Consider $n$ independent Brownian motions with transition probability
\eqref{transitionprob}--\eqref{variance} conditioned so that
\begin{itemize}
\item the paths are non-intersecting in the time interval $(0,1)$,
\item $n_{k,l}$ of the paths start at $a_k$ and end at $b_l$,
for each $k = 1, \ldots, p$ and $l =1, \ldots, q$.
\end{itemize}
Assume that as $n \to \infty$, the finite $n$ transition numbers $n_{k,l}/n$
converge to $t_{k,l}$, and that the corresponding graph $G$ is connected.

Let $t \in (0,1)$. Then there exists a $T^* = T^*(t) > 0$ so that for all $T
\in (0, T^*)$ the limiting mean distribution of the positions of the paths at
time $t$ exists, and is supported on the union of $p+q-1$ disjoint intervals
$\bigcup_{i=1}^{p+q-1} [\alpha_i, \beta_i]$ with a density $\rho_i$ on the
$i$th interval.

The vector of measures $(\mu_1, \ldots, \mu_{p+q-1})$ where $d\mu_i(x) =
\rho_i(x) dx$ for $i=1, \ldots, p+q-1$ is the minimizer of a vector equilibrium
problem that will be described in the next subsection.
\end{theorem}

The proof of Theorem \ref{theorem:limitingdistribution} will be based on a
Deift-Zhou steepest descent analysis of the Riemann-Hilbert problem in Section
\ref{subsection:RHintro}. The details of the steepest descent analysis will be
described in Sections
\ref{section:normalizationatinfinity}--\ref{section:remainingstepssteepestdescent},
and the proof of Theorem~\ref{theorem:limitingdistribution} will then be given
in Section~\ref{section:proofmaintheorems}.

\begin{remark} \label{remark:bleherkuijlaars}
The special case where $p=1$, $q=2$ and $m_1=m_2=n/2$ was studied
in \cite{ABK,BK2,BK3}. In \cite{BK2} the energy functional
\eqref{energyminimization:mixedangelesco}  is given for this
special case, but it was not used in the further analysis. Instead
the results in \cite{BK2} were stated and proved in terms of an
algebraic curve (Pastur's equation).
\end{remark}

\begin{remark} \label{remark:generalexternalfields}
As already mentioned, for the case $p=1$ and $q \geq 2$
the non-intersecting Brownian motions are distributed like the
eigenvalues of the Gaussian unitary ensemble with external source.

There are more general random matrix ensembles with external source,
which however do not have an equivalent interpretation in terms of non-intersecting paths.
The case of a quartic potential was studied in \cite{Mcl}.
An analogue of Theorem \ref{theorem:limitingdistribution} is valid
in this case as well, provided that a suitable analogue of the temperature
$T$ is sufficiently small.
\end{remark}

Note that Theorem \ref{theorem:limitingdistribution} gives only a result for
$T$ sufficiently small, but it does not specify how small $T$ should be. We
expect that the theorem remains valid for all temperatures $T$ that are such
that at time $t$ we have the maximal number of groups of paths. At a critical
temperature $T_{\crit}(t)$ we expect that two (or maybe more) neighboring
intervals $[\alpha_i,\beta_i]$ and $[\alpha_{i+1},\beta_{i+1}]$ merge and a Pearcey phase transition occurs.
However we were unable to prove this.

\subsection{Equilibrium problem}
\label{subsection:generalstatementmainresults}

The logarithmic energy of a measure $\mu$ on $\mathbb R$ is defined as usual
by
\begin{equation}\label{deflogarithmicenergy}
I(\mu) = \iint \log\frac{1}{|x-y|} \, d\mu(x)\, d\mu(y).
\end{equation}
The mutual energy of two measures $\mu$, $\nu$ is defined by
\begin{equation}\label{defmutualenergy}
I(\mu,\nu) = \iint \log\frac{1}{|x-y|} \, d\mu(x)\, d\nu(y).
\end{equation}

We write for $i = 1,\ldots, p+q-1$,
\begin{equation} \label{defxit}
    x_i(t) = (1-t) a_{k(i)} + t b_{l(i)}, \qquad 0 < t < 1
    \end{equation}
and given $t \in (0,1)$ we define the quadratic functions
\begin{equation}\label{defexternalfield}
V_i(x) = \frac{1}{2 t(1-t)} (x-x_i(t))^2, \qquad x \in \mathbb R,
\end{equation}
for $i=1, \ldots, p+q-1$. These functions will play the role of external fields
in the vector equilibrium problem that is relevant for our problem.

\begin{definition} \label{def:equilibriumproblem}
Fix $t\in (0,1)$ and let $T > 0$. Consider
the energy functional
\begin{equation}\label{energyminimization:mixedangelesco}
E(\mu_1,\ldots,\mu_{p+q-1}) = \sum_{i,j=1}^{p+q-1}
    a_{i,j}I(\mu_i,\mu_j)+  \frac{1}{T} \sum_{i=1}^{p+q-1} \int V_i(x) \, d\mu_i(x),
\end{equation}
where $V_i$ is defined in \eqref{defexternalfield},
and the interaction matrix $A = (a_{i,j})$ has entries
\begin{equation}\label{def:interactionmatrix}
a_{i,j} = \left\{\begin{array}{ll}
    1 & \textrm{ if } i=j, \\[5pt]
    \frac{1}{2} & \textrm{ if } i\neq j \textrm{ and } k(i)=k(j) \textrm{ or } l(i)=l(j), \\[5pt]
    0 & \textrm{ otherwise.}
\end{array}\right.
\end{equation}
Then the vector equilibrium problem consists in minimizing the energy functional
\eqref{energyminimization:mixedangelesco} over all vectors of positive measures
$(\mu_1,\ldots,\mu_{p+q-1})$ supported on the real line for which
\begin{equation}\label{totalmassmui}
\int d\mu_i = t_{k(i),l(i)}, \qquad \text{for } i=1,\ldots,p+q-1.
\end{equation}
\end{definition}

One may understand the energy minimization problem in Definition
\ref{def:equilibriumproblem} as follows. To each of the edges of the graph $G = (V,E,t)$
we associate a measure $\mu_i$, $i=1,\ldots,p+q-1$, of total mass equal
to the weight $t_{k(i),l(i)}$ of that edge. This measure represents a distribution of
charged particles on the real line that repel each other due
to the diagonal term
\[ a_{i,i} I(\mu_i, \mu_i) = I(\mu_i, \mu_i) \]
in \eqref{energyminimization:mixedangelesco}.
For the particles of different measures $\mu_i$, $\mu_j$, $i\neq j$, there are two possibilities.
The first case is when the $(i,j)$ entry of \eqref{def:interactionmatrix} equals $1/2$.
This happens if the edges corresponding to $i$ and $j$ are adjacent in the graph $G$.
Then there is repulsion between the measures $\mu_i$ and $\mu_j$
but with a strength that is only half as strong as the repulsion for each individual measure.
The second case is when the $(i,j)$ entry of \eqref{def:interactionmatrix} equals zero. In
that case there is no direct interaction between the measures $\mu_i$ and $\mu_j$.

The last term of \eqref{energyminimization:mixedangelesco} is a sum
of external field terms due to the action of the external field $\frac{1}{T} V_i(x)$
on the measure $\mu_i$. The energy minimizer
$(\mu_1,\ldots,\mu_{p+q-1})$ in Definition~\ref{def:equilibriumproblem} then
corresponds to the equilibrium distribution of charged particles under the
energy functional \eqref{energyminimization:mixedangelesco}.

\begin{proposition} \label{prop:positivedefinite}
The interaction matrix $A$ is positive definite.
\end{proposition}

\begin{proof}
It is easy to check that the interaction matrix $A$ is
equal to
\begin{equation}     \label{factorization:edgeconnection}
    A = \frac{1}{2} B^T \, B
    \end{equation}
where $B$ is the incidence matrix of the graph $G$. That is,
we choose a numbering $k=1, \ldots, p+q$ of the vertices, and then we have
\[ B = (b_{k,i})_{k=1, \ldots, p+q, \, i = 1, \ldots, p+q-1} \]
where $b_{k,i} = 1$ if vertex $k$ is incident to edge $i$,
and $0$ otherwise.

From \eqref{factorization:edgeconnection} we get that $A$ is
positive semi-definite, and for any column vector $\vec{x}$ of length $p+q-1$
we have
\begin{equation} \label{positivesemidef}
 \vec{x}^T  A \vec{x} = \frac{1}{2} \| B\vec{x} \|^2 \geq 0.
\end{equation}

Now assume that $B\vec{x} = \vec{0}$. Consider a leaf of $G$, i.e., a vertex
which is incident to exactly one edge. Then $B$ has exactly one zero in the row
corresponding to this vertex, and from $B\vec{x} = \vec{0}$ it follows that the
component of $\vec{x}$ corresponding to the edge that is incident to the leaf
vanishes. Since $G$ is a tree (see Proposition \ref{lemma:forest}(c)) we can
then gradually undress $G$ by peeling off leafs one by one and we conclude in
this way that all components of $\vec{x}$ are equal to $0$. Thus
$\vec{x}=\vec{0}$ if $B \vec{x} = \vec{0}$, which implies in view of
\eqref{positivesemidef} that $A$ is positive definite.
\end{proof}

\begin{corollary}
The vector equilibrium problem of Definition \ref{def:equilibriumproblem}
has a unique solution $(\mu_1, \ldots, \mu_{p+q-1})$ and each measure $\mu_i$
is compactly supported.
\end{corollary}

\begin{proof}
 The interaction matrix is positive definite by Proposition
\ref{prop:positivedefinite}. The external fields $V_i$ in
the energy functional \eqref{energyminimization:mixedangelesco} have
enough increase at $\pm \infty$ so that standard arguments
of potential theory as in \cite{Dei,NS,SaffTotik} can be used
to establish the existence and uniqueness of the minimizer as well as the
fact that each measure $\mu_i$ is supported on a compact set.
\end{proof}

Our next theorem describes the structure of the solution of the vector
equilibrium problem for small $T$.

\begin{theorem}\label{theorem:energyminimizerexists}
Fix $t\in (0,1)$, and let $(t_{k,l})$ be a matrix of transition numbers. Then
there exists $T^{*} > 0$ (the same $T^{*}$ that makes Theorem
\ref{theorem:limitingdistribution} work) so that for every $T \in (0,T^{*})$
the following holds.
\begin{itemize}
\item[\rm (a)] Each $\mu_i$ is supported on an interval
\[ \supp(\mu_i) = [\alpha_i,\beta_i], \qquad i=1,\ldots,p+q-1. \]
The intervals $[\alpha_i,\beta_i]$ are pairwise disjoint and satisfy
\begin{equation}\label{orderoftheintervals}
    \beta_{i+1} < \alpha_{i}, \qquad i=1,\ldots,p+q-2.
    \end{equation}
\item[\rm (b)] The measure $\mu_i$ has a density $\rho_i$ with respect to Lebesgue measure
which is real analytic and positive in the open interval $(\alpha_i, \beta_i)$
and vanishes like a square root at the endpoints of $[\alpha_i,\beta_i]$,
i.e., there exist non-zero constants $\rho_i^{(1)}$ and
$\rho_i^{(2)}$ such that
\begin{align}\label{squarerootbehavior:alpha}
    \rho_i(x) = \rho_i^{(1)}\sqrt{x-\alpha_i}+O((x-\alpha_i)^{3/2}) & \quad\textrm{ as } x\downarrow\alpha_i,\\
    \label{squarerootbehavior:beta}
    \rho_i(x) = \rho_i^{(2)}\sqrt{\beta_i-x}+O((\beta_i-x)^{3/2}) & \quad \textrm{ as }
    x\uparrow\beta_i.
\end{align}
\end{itemize}
\end{theorem}

The proof of Theorem \ref{theorem:energyminimizerexists} will be given in
Section \ref{subsection:proofequilibriumproblem}.

\subsection{Special cases}

In Subsections \ref{subsection:angelesco}--\ref{subsection:mixedangelesco}, our
main Theorem \ref{theorem:limitingdistribution} will be illustrated for some
special cases.

\subsubsection{The case $p=1$: Angelesco-type interaction}
\label{subsection:angelesco}

\begin{figure}[htbp]\vspace{-8mm}
\begin{center}
        \subfigure{\label{blablabla}}\includegraphics[scale=0.5]{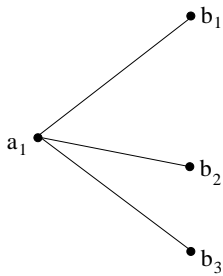}
\end{center}\vspace{-1mm}
\caption{A graph $G$ with $p=1$ starting points and $q=3$ ending points.}
\label{fig:graphAngelesco}
\end{figure}

For the case $p=1$ of one starting point, and an arbitrary number $q$ of ending
points, the graph $G$ has a single vertex $a_1$ on the left which is connected
to each of the vertices $b_1, \ldots, b_q$ on the right. For example, if $q=3$
the graph has the form shown in Figure \ref{fig:graphAngelesco}.

The energy functional
\eqref{energyminimization:mixedangelesco}--\eqref{def:interactionmatrix}
is then equal to
\begin{equation}\label{energyminimization:angelesco}
E(\mu_1,\ldots,\mu_{p+q-1}) := \sum_{i=1}^{p+q-1} I(\mu_i) +
\frac{1}{2}\sum_{i\neq j} I(\mu_i,\mu_j) + \frac{1}{T} \sum_{i=1}^{p+q-1} \int  V_i(x)\, d\mu_i(x).
\end{equation}
The functional \eqref{energyminimization:angelesco} is exactly the one familiar
from the theory of Angelesco systems \cite{GRS}. All off-diagonal entries in the
interaction matrix $A$ in \eqref{def:interactionmatrix} are equal to $1/2$.
For the example in
Figure \ref{fig:graphAngelesco} the interaction matrix is
\begin{equation}\label{interactionmatrix:angelesco}
A =
\begin{pmatrix}
1 & 1/2 & 1/2\\
1/2 & 1 & 1/2\\
1/2 & 1/2 & 1\\
\end{pmatrix}.
\end{equation}

\subsubsection{The \lq zigzag\rq\ case: nearest neighbor interaction}
\label{subsection:zigzag}

\begin{figure}[htbp]
\begin{center}\vspace{-8mm}
        \subfigure[]{\label{blablabla1}}\includegraphics[scale=0.5]{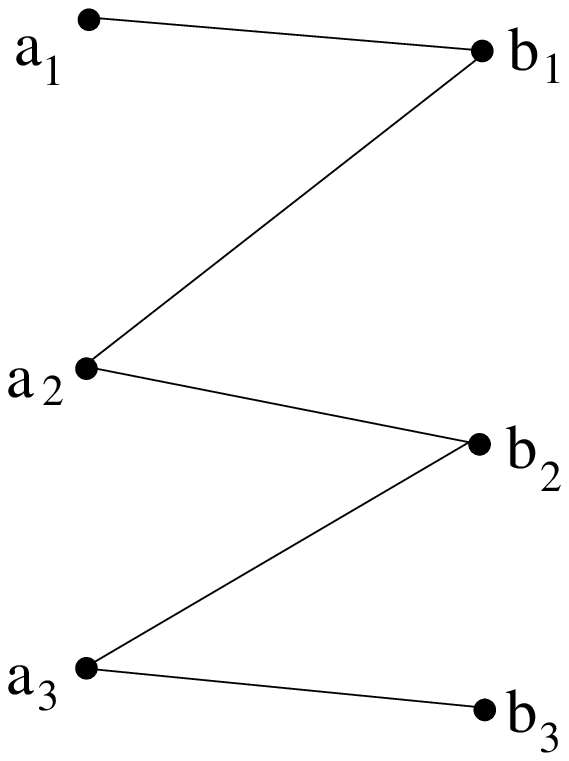}\hspace{10mm}
        \subfigure[]{\label{blablabla2}}\includegraphics[scale=0.5]{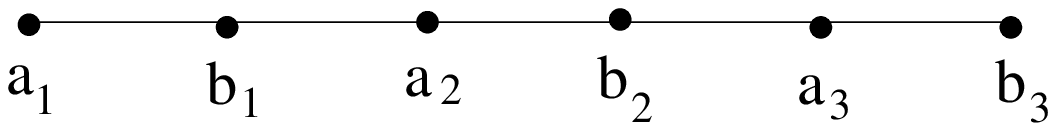}
\end{center}\vspace{-1mm}
\caption{For $p=q=3$, Figure \ref{blablabla1} shows a graph $G$ which has
zigzag form. Figure \ref{blablabla2} shows the same graph written as a chain.}
\label{fig:graphzigzag}
\end{figure}

Next we consider the case where $p=q$ and the corresponding lattice path
follows a zigzag line. The graph $G$ is then just a chain of
vertices: see Figure \ref{fig:graphzigzag}.

The energy functional \eqref{energyminimization:mixedangelesco} now takes the
form
\begin{equation}\label{energyminimization:tridiagonal}
E(\mu_1,\ldots,\mu_{p+q-1}) := \sum_{i=1}^{p+q-1} I(\mu_i) + \sum_{i=1}^{p+q-2}
I(\mu_i,\mu_{i+1}) + \frac{1}{T} \sum_{i=1}^{p+q-1} \int V_i(x) \, d\mu_i(x).
\end{equation}
We see that the interaction matrix of \eqref{energyminimization:tridiagonal} is
tridiagonal with diagonal entries equal to $1$, and entries on the first sub-
and superdiagonal equal to $1/2$. For the example in Figure
\ref{fig:graphzigzag} the interaction matrix equals
\begin{equation}\label{interactionmatrix:tridiagonal}
A =
\begin{pmatrix}
1&1/2&0&0&0\\
1/2&1&1/2&0&0\\
0&1/2&1&1/2&0\\
0&0&1/2&1&1/2\\
0&0&0&1/2&1\\
\end{pmatrix}.
\end{equation}
Note that the tridiagonal structure of the interaction matrix means that there
is only nearest neighbor interaction. The neighboring measures repel
each other, since all signs in the interaction matrix are positive.

In a Nikishin system the interaction matrix is also tridiagonal, but the
off-diagonal entries are $-1/2$ instead of $1/2$, see \cite{GRS}.

\subsubsection{The general case: block nearest neighbor interaction with
Angelesco-type blocks} \label{subsection:mixedangelesco}

Finally we consider the graph $G$ in Figure \ref{fig:graph2}. In this case the
interaction matrix equals
\begin{equation}\label{interactionmatrix:generalexample}
A = \begin{pmatrix}
1&1/2&0&0&0\\
1/2&1&1/2&0&0\\
0&1/2&1&1/2&1/2\\
0&0&1/2&1&1/2\\
0&0&1/2&1/2&1
\end{pmatrix}.
\end{equation}
Note that this interaction matrix is a mixture of the nearest neighbor and
Angelesco interaction matrices. More precisely one could say that
\eqref{interactionmatrix:generalexample} has \lq block\rq\ nearest neighbor
interaction, where each of the blocks in turn has an Angelesco-type
interaction, and with subsequent blocks intersecting in exactly one entry. For
the matrix in \eqref{interactionmatrix:generalexample} these building blocks
are
$$ \begin{pmatrix}1 & 1/2\\1/2& 1\end{pmatrix},\quad \begin{pmatrix}1 & 1/2\\1/2& 1\end{pmatrix},\quad
\begin{pmatrix}1 & 1/2 & 1/2\\1/2& 1 & 1/2\\ 1/2&1/2&1\end{pmatrix}.
$$

\subsubsection{Comparison with \cite{GRS}}

A graph-based vector equilibrium problem was also considered in the work of
Gonchar, Rakhmanov, and Sorokin \cite{GRS}. The rule for building the
interaction matrix from a graph $G$ is similar to the one in this paper. The
vector equilibrium problem is also labeled by the edges of a graph, which in
\cite{GRS}, however, is a directed rooted tree. The off-diagonal entries of the
interaction matrix are non-zero if the corresponding two edges have a common
vertex. The entry is $1/2$ if the two edges have a common initial vertex, and
$-1/2$ if the final vertex of one edge agrees with the initial vector of the
other.

In our case, we only have $1/2$ since the latter situation
cannot happen.

\subsection{About the proof of Theorem \ref{theorem:limitingdistribution}}
\label{subsection:RHintro}

The proof of Theorem \ref{theorem:limitingdistribution} uses the connection of
non-intersecting Brownian motions with prescribed starting and ending points
with a determinantal point process and an associated Riemann-Hilbert problem.
We recall this connection.

\subsubsection{Determinantal point process}

The positions at time $t \in (0,1)$ of $n$ non-intersecting Brownian motions,
starting at distinct $a_j$, $j=1, \ldots, n$, and ending at distinct positions
$b_j$, $j=1, \ldots, n$, have the joint probability density function
\begin{equation*} \label{jointPDF}
    \frac{1}{Z} \det \left( P(t,a_i, x_j) \right)_{i,j=1}^n
    \cdot \det \left( P(1-t, x_i, b_j) \right)_{i,j=1}^n,
    \end{equation*}
with the transition probability density $P$ defined in \eqref{transitionprob}
and with $Z$ a normalization constant. This is a consequence of a theorem of
Karlin and McGregor \cite{KMcG}. (For applications of the discrete version of
the Karlin-McGregor theorem see e.g.\ \cite{Joh1}.) In the confluent limit
where $n_k$ of the starting positions come together at $a_k$, $k=1,\ldots,p$,
and $m_l$ of the ending positions come together at $b_l$ for $l=1,\ldots,q$,
the joint p.d.f.\ for the positions of the paths at time $t$ can be written as
\begin{equation} \label{jointPDFcoalesce}
    \mathcal P(x_1, \ldots, x_n) =
    \frac{1}{\tilde{Z}} \det \left( f_i(x_j) \right)_{i,j=1}^n
        \cdot \det \left( g_i(x_j) \right)_{i,j=1}^n,
\end{equation}
for certain functions $f_i$, $g_i$, that are built out
of the $p+q$ functions
\begin{align}
\label{gaussianw1}
    w_{1,k}(x) & = e^{-\frac{1}{2t\sigma^2}(x-a_k)^2},\quad k=1,\ldots,p,\\
\label{gaussianw2}
    w_{2,l}(x) & = e^{-\frac{1}{2(1-t) \sigma^2}(x-b_l)^2},\quad
    l=1,\ldots,q,
\end{align}
see e.g.\ \cite{DK2}.

The p.d.f.\ \eqref{jointPDFcoalesce} defines a determinantal point process (in
fact a biorthogonal ensemble, see \cite{Bor}) with a correlation kernel $K(x,y)$
that is such that
\begin{equation} \label{pdfandkernel1}
    \mathcal P(x_1, \ldots, x_n)
            = \frac{1}{n!} \det \left( K(x_i,x_j) \right)_{i,j=1}^n
\end{equation}
and for each $m = 1, \ldots, n$,
\begin{equation} \label{pdfandkernel2}
    \int \cdots \int \mathcal P(x_1, \ldots, x_n) \, dx_{m+1} \cdots d x_n
    = \frac{(n-m)!}{n!} \det \left( K(x_i,x_j) \right)_{i,j=1}^m.
\end{equation}
In particular for $m=1$, we have that
\[ \frac{1}{n} K(x,x) \]
is the mean density of paths.

\subsubsection{Riemann-Hilbert problem}

The kernel $K$ can be
described in terms of the following Riemann-Hilbert problem (RH problem)
introduced in \cite{DK2}.

\begin{rhp}\label{RHP:original}
The RH problem consists in finding a matrix-valued function
$Y : \mathbb C \setminus \mathbb R \to \mathbb C^{(p+q)\times(p+q)}$ such that
\begin{itemize}
\item[\rm (1)] $Y$ is analytic in $\mathbb C \setminus\mathbb R $;
\item[\rm (2)] For $x \in \mathbb R$, the limiting values
\[ Y_+(x) = \lim_{z \to x, \,  \Im z > 0} Y(z), \qquad
    Y_-(x) = \lim_{z \to x, \, \Im z < 0} Y(z) \]
exist and satisfy
\begin{equation}\label{jumpsY10}
Y_{+}(x) = Y_{-}(x)
\begin{pmatrix} I_p & W(x)\\
0 & I_q
\end{pmatrix},
\end{equation}
where $I_k$ denotes the identity matrix of size $k$,  and where $W(x)$ denotes the
rank-one matrix (outer product of two vectors)
\begin{equation}\label{defWblock0}
W(x) = \begin{pmatrix} w_{1,1}(x) \\ \vdots \\ w_{1,p}(x)
\end{pmatrix}\begin{pmatrix} w_{2,1}(x) & \cdots & w_{2,q}(x)
\end{pmatrix}
\end{equation}
with $w_{1,k}(x)$, $k=1, \ldots,p$, and  $w_{2,l}(x)$, $l=1, \ldots, q$ given
by \eqref{gaussianw1} and \eqref{gaussianw2}.
\item[\rm (3)] As $z\to\infty$, we have that
\begin{equation}\label{asymptoticconditionY0} Y(z) =
    (I_{p+q}+O(1/z))\diag\left(z^{n_1},\ldots,z^{n_p},z^{-m_1},\ldots,z^{-m_q}\right).
\end{equation}
\end{itemize}
\end{rhp}

The RH problem has a unique solution that can be described in terms of multiple
Hermite polynomials. This is shown in \cite{DK2}, generalizing results in
\cite{FIK,VAGK}.
According to \cite{DK2} the correlation kernel  is  expressed in terms of the
solution to the RH problem as
\begin{equation} \label{correlationkernel}
    K(x,y) = \frac{1}{2\pi i(x-y)}\begin{pmatrix} 0 & \cdots & 0 & w_{2,1}(y) &
    \ldots & w_{2,q}(y)\end{pmatrix} Y_{+}^{-1}(y)Y_{+}(x)\begin{pmatrix}
    w_{1,1}(x)\\ \vdots \\
    w_{1,p}(x)\\ 0 \\ \vdots \\ 0 \end{pmatrix}.
\end{equation}
It is also worth noticing that if $Y_{1,1}(z)$ denotes the top leftmost $p$ by
$p$ block of the RH matrix $Y(z)$, then the determinant of $Y_{1,1}(z)$ equals
the average characteristic polynomial
$$\det Y_{1,1}(z) = \mathbb E\left[\prod_{j=1}^{n}(z-x_j)\right],$$ where the expectation $\mathbb E$ is taken
according to the joint probability density \eqref{jointPDFcoalesce}. This was
shown in \cite{BK} for the case $p=1$ and in \cite{Del} for the case of general
$p$ and $q$.

\subsubsection{Asymptotic analysis}

We analyze the RH problem \ref{RHP:original} in the limit described in Theorem
\ref{theorem:limitingdistribution}. That is, we take $\sigma^2 = T/n$ and we
let $n_k \to \infty$, $m_l \to \infty$, so that
\[ \frac{n_{k,l}}{n} \to t_{k,l} \qquad \text{as } n \to \infty. \]

If, for each $n$, we denote the correlation kernel \eqref{correlationkernel}
by $K_n$, then the limiting density of paths
is
\begin{equation} \label{correlationkernellimit}
    \rho(x) = \lim_{n\to\infty} \frac{1}{n} K_n(x,x), \qquad x \in \mathbb R.
\end{equation}

The proof of Theorem \ref{theorem:limitingdistribution} will be based on a
steepest descent analysis of the RH problem~\ref{RHP:original}. As a byproduct
of this analysis we can also show that the local scaling limits of the limiting
distribution of the Brownian motions are those familiar from random matrix
theory, i.e., they are described in terms of the sine kernel in the bulk and
the Airy kernel at the edge. We will not discuss this any further and refer to
the papers \cite{BK2,BK3,DKV,Ora} for a more detailed analysis in a similar
context.

\subsection{Outline of the rest of the paper}

The rest of this paper is organized as follows. In Section
\ref{section:propertiesequilibriumproblem} we establish general properties of
the equilibrium problem, in particular leading to the proof of
Theorem~\ref{theorem:energyminimizerexists}. In Section
\ref{section:xifunctions} we define the $\xi$-functions, $\lam$-functions and
the associated Riemann surface. These functions are used in Section
\ref{section:normalizationatinfinity} to normalize the RH problem at infinity.
In Section~\ref{section:globallenses} we apply Gaussian elimination to the RH
problem by opening global lenses, thereby making the RH problem locally of size
$2\times 2$. This construction is fully systematic and may have interest in its
own right. The remaining steps of the RH steepest descent analysis are
described in Section \ref{section:remainingstepssteepestdescent}. Finally, in
Section \ref{section:proofmaintheorems} we prove the main
Theorem~\ref{theorem:limitingdistribution}.

\section{Properties of the equilibrium measures}
\label{section:propertiesequilibriumproblem}

In this section we study the vector equilibrium problem of Definition
\ref{def:equilibriumproblem}. In particular we prove Theorem
\ref{theorem:energyminimizerexists}.

\subsection{Proof of Theorem \ref{theorem:energyminimizerexists}}
\label{subsection:proofequilibriumproblem}

The main difficulty is in the proof of part (a)
of Theorem \ref{theorem:energyminimizerexists}.
For this we use the following lemma. Recall that $t \in (0,1)$
is fixed.

\begin{lemma} \label{lemma:energyminimizer}
For every $\varepsilon > 0$ and $\tau > 0$, there exists $T_{\varepsilon} > 0$
so that for every $T \in (0,T_{\varepsilon})$ the following holds.

If $(\mu_1, \ldots, \mu_{p+q-1})$ is the minimizer for the energy
functional \eqref{energyminimization:mixedangelesco} under
the normalization \eqref{totalmassmui} with
\begin{equation} \label{mintkili}
    \min_i t_{k(i),l(i)} \geq \tau > 0,
    \end{equation}
then the support of $\mu_i$
is contained in $[x_i(t)-\varepsilon, x_i(t) + \varepsilon]$
for every $i = 1, \ldots, p+q-1$.
\end{lemma}

\begin{proof}
We are going to prove the following upper and lower bounds for
the quantity
\[ E^* := E(\mu_1, \ldots, \mu_{p+q-1}) \]
(which depends on the chosen normalization \eqref{totalmassmui}).
\begin{enumerate}
\item[(a)] There exist constants $C_1 > 0$ and $C_{2,\delta}$ such
that for every $\delta > 0$ and $T > 0$ we have
\begin{equation} \label{upperboundonE}
    E^* \leq C_{2,\delta} + C_1 \frac{\delta^2}{T}.
    \end{equation}
The constant $C_1$ is independent of $\delta$ and $T$, and
$C_{2,\delta}$ is independent of $T$.
\item[(b)] There exist positive constants $C_3$, and $C_4$, such
that for every $T \in (0,1]$ and every $i = 1, \ldots, p+q-1$,
we have
\begin{equation} \label{lowerboundonE}
    E^* \geq - C_3 + \frac{C_4}{T} (x-x_i(t))^2,
        \qquad \text{for } x \in \supp(\mu_i).
    \end{equation}
\end{enumerate}
In addition, all constants $C_1$, $C_{2,\delta}$, $C_3$ and $C_{4}$
can be taken independently of the normalization \eqref{totalmassmui},
and only $C_4$ depends on $\tau$.

\paragraph{Proof of (a):} We may assume that $\delta > 0$ is  sufficiently small so that
\[ \min_i \left(x_i(t) - x_{i+1}(t)\right) \geq 2 \delta. \]
Then the intervals $[x_i(t)-\delta, x_i(t) +\delta]$ are mutually disjoint.

We fix transition numbers $t_{k(i),l(i)} \geq 0$ so that $\sum_i t_{k(i),l(i)} = 1$.
On each of the intervals $[x_i(t)-\delta, x_i(t) + \delta]$ we choose a measure
$\lambda_i$ of finite energy with total mass $1$
(for example, an appropriately rescaled and centered
semi-circle law will do) and we put $\nu_i = t_{k(i),l(i)} \lambda_i$.
We choose the measures $\nu_i$ independent of $T$ and so
\begin{equation} \label{upperboundonE20}
    \sum_{i,j=1}^{p+q-1} a_{i,j} I(\nu_i,\nu_j)
    \end{equation}
does not depend on $T$. The sum \eqref{upperboundonE20} depends
continuously on the numbers $t_{k(i),l(i)}$ and so
\begin{equation} \label{upperboundonE2}
    C_{2,\delta} = \max \sum_{i,j=1}^{p+q-1} a_{i,j} I(\nu_i, \nu_j)
\end{equation}
exists, where the maximum is taken over all choices $t_{k(i),l(i)} \geq 0$
with $\sum_i t_{k(i),l(i)} = 1$.

Since $V_i(x) \leq \frac{\delta^2}{2t(1-t)}$ on $[x_i(t)-\delta, x_i(t) + \delta]$,
we also have
\[ \int V_i(x) d\nu_i(x) \leq \frac{\delta^2}{2t(1-t)} \]
and
\begin{equation} \label{upperboundonE3}
    \sum_{i=1}^{p+q-1} \int V_i(x) \, d\nu_i(x) \leq  C_1 \delta^2
    \end{equation}
with a constant $C_1$ that is independent of $T$ and $\delta$, and also of
the $t_{k(i),l(i)}$.
Then part (a) follows from \eqref{upperboundonE2}, \eqref{upperboundonE3},
and the fact that
\[  E^* \leq E(\nu_1, \ldots, \nu_{p+q-1})
    = \sum_{i,j=1}^{p+q-1} a_{i,j} I(\nu_i,\nu_j) + \frac{1}{T} \sum_{i=1}^{p+q-1} \int V_i(x) d\nu_i(x). \]

\paragraph{Proof of (b):}

For each $i$ we take a measure $\nu_i$ of total mass $\| \nu_i \| = \| \mu_i \|$
and satisfying for some constant $K \geq 1$,
\begin{equation} \label{lowerboundonE2}
    \nu_i \leq K \mu_i, \qquad \text{for } i = 1, \ldots, p+q-1.
    \end{equation}
Then $(1-\lambda) \mu_i + \lambda \nu_i$ is a positive measure
for every $\lambda \in [-1/K, 1]$, and so
\begin{equation} \label{lowerboundonE3}
    E^* \leq E\left((1-\lambda)\mu_1+\lambda\nu_1,\ldots,(1-\lambda)\mu_{p+q-1}+\lambda\nu_{p+q-1}\right)
\end{equation}
for every $\lambda \in [-1/K,1]$, with equality for $\lambda = 0$.
Thus the $\lambda$-derivative of the right-hand side
of \eqref{lowerboundonE3} vanishes for $\lambda = 0$, which leads to
\begin{equation}\label{lowerboundonE4}
    E^* =
    \sum_{i,j=1}^{p+q-1} a_{i,j}I(\mu_i,\nu_j) +\frac{1}{2T} \sum_{i=1}^{p+q-1} \int V_i(x) \,
    \left( d\mu_i(x) + d\nu_i(x) \right).
    \end{equation}

From the elementary inequality $|x-y| \leq \sqrt{x^2+1} \sqrt{y^2+1}$ it follows that
\begin{align*}
    I(\mu_i, \nu_j) & = \iint \log \frac{1}{|x-y|} d\mu_i(x) d\nu_j(y) \\
    & \geq - \frac{1}{2} \iint \log(x^2 + 1) d\mu_i(x) d\nu_j(y)
         - \frac{1}{2} \iint \log(y^2+1) d\mu_i(x) d\nu_j(y) \\
    & \geq - \frac{1}{2} \int \log(x^2+1) \left(d\mu_i(x) + d\nu_j(x) \right)
    \end{align*}
where for the last inequality we used the facts that $\log(x^2+1) \geq 0$ and $\| \nu_j\| \leq 1$, $\|\mu_i \|\leq 1$.
Since all $a_{i,j} \geq 0$ it then follows from
\eqref{lowerboundonE4} that
\begin{align} \nonumber
    E^* & \geq - \frac{1}{2}
        \sum_{i,j=1}^{p+q-1} a_{i,j}  \int \log(x^2+1) \left(d\mu_i(x) + d\nu_j(x)\right)
            + \frac{1}{2T} \sum_{i=1}^{p+q-1} \int V_i(x) \left(d\mu_i(x) + d\nu_i(x) \right) \\
    & = \label{lowerboundonE5}
    \frac{1}{2T} \sum_{i=1}^{p+q-1} \int \left( V_i(x) -  T \sum_{j=1}^{p+q-1} a_{i,j} \log(x^2+1)\right)\, (d\mu_i(x)+d\nu_i(x)).
\end{align}

For $i=1, \ldots, p+q-1$ and $T \leq 1$, the functions
$V_i(x) - T \sum_{j=1}^{p+q-1} a_{i,j} \log(x^2+1)$ are bounded from below
and assume their minimum value in a fixed compact interval (independent of $i$ and $T \leq 1$).
This is clear from the explicit form \eqref{defexternalfield} of the functions $V_i$,
and from the fact that $\sum_{j=1}^{p+q-1} a_{i,j} \leq \frac{p+q}{2}$
for each $i$. Thus
\[ V_i(x) - T \sum_{j=1}^{p+q-1} a_{i,j} \log (x^2+1) \geq \frac{1}{2} V_i(x) -C_0 T
    \geq - C_0 T, \qquad x \in \mathbb R, \]
for some constant $C_0 > 0$ independent of $T \leq 1$ and $i$.

Using this in \eqref{lowerboundonE5} we find that for each $i =1, \ldots, p+q-1$,
\begin{equation} \label{lowerboundonE6}
    E^* \geq
    -C_3 + \frac{1}{4T} \int V_i(x) \, d\nu_i(x),
    \end{equation}
where $C_3 > 0$ is a constant independent of $i$ and $T \leq 1$.
The inequality \eqref{lowerboundonE6} (with the same constant $C_3$)
holds for all measures $\nu_i$ with $\| \nu_i \| = t_{k(i),l(i)}$
and satisfying \eqref{lowerboundonE2} for some $K$.
For every $x$ in the support of $\mu_i$, we can approximate
the point mass $\|\mu_i\| \, \delta_x$ by such $\nu_i$.
It thus follows from \eqref{lowerboundonE6} that
\[ E^* \geq -C_3 + \frac{1}{4T} V_i(x) \, t_{k(i),l(i)},
    \qquad \textrm{for } x \in \supp(\mu_i). \]
If  $t_{k(i),l(i)} \geq \tau$ as in \eqref{mintkili},
then we  obtain \eqref{lowerboundonE} with
 \[ C_4 = \frac{1}{4} \cdot \frac{1}{2t(1-t)} \cdot \tau. \]

\paragraph{Conclusion of the proof:}

Combining \eqref{upperboundonE} and \eqref{lowerboundonE}
we find that there exist positive constants $C_5 = C_1/C_4$ and $C_{6,\delta} = (C_{2,\delta} + C_3)/C_4$ so that
\[ (x-x_i(t))^2 \leq C_5 \delta^2 + C_{6,\delta} T, \qquad
    \text{for } x\in \supp(\mu_i), \]
 for every $i =1, \ldots, p+q-1$, $T \leq 1$ and $\delta > 0$.

Let $\varepsilon > 0$ be given. Choose first $\delta > 0$ so that $C_5 \delta^2
\leq \frac{1}{2} \varepsilon^2$ and then choose $T_{\varepsilon} \in (0,1]$ so
that $C_{6,\delta} T_{\varepsilon} \leq \frac 12\varepsilon^2$. Then for every
$i = 1, \ldots, p+q-1$ and $T \leq T_{\varepsilon}$,
\[ (x-x_i(t))^2 \leq \varepsilon^2, \qquad
    \text{for } x\in \supp(\mu_i), \]
and  Lemma \ref{lemma:energyminimizer} follows.
\end{proof}

Having Lemma \ref{lemma:energyminimizer}
the proof of Theorem \ref{theorem:energyminimizerexists} is
rather straightforward.

\begin{proof} \textbf{(Proof of Theorem \ref{theorem:energyminimizerexists})}

(a) Let $\varepsilon > 0$ be such that the intervals
$[x_i(t) - \varepsilon, x_i(t) + \varepsilon]$, $i=1, \ldots, p+q-1$, are disjoint.
From Lemma \ref{lemma:energyminimizer} we know that there exists $T_{\epsilon} > 0$
so that for $T < T_{\varepsilon}$ the support of $\mu_i$ is contained
in $[x_i(t) - \varepsilon, x_i(t) + \varepsilon]$.

Take $i =1, \ldots, p+q-1$, and fix the other measures $\mu_j$, $j\neq i$.
From \eqref{defmutualenergy}, \eqref{energyminimization:mixedangelesco}
we see that the measure $\mu_i$ is the equilibrium measure in an effective external field
\begin{equation}\label{netexternalfield}
    \sum_{j\neq i}a_{i,j}\int \log\frac{1}{|x-y|}\, d\mu_j(y) +
    \frac{1}{T} V_i(x),
\end{equation}
and it is also the minimizer if we restrict to measures with total mass $t_{k(i),l(i)}$
that are supported on $[x_i(t)-\varepsilon, x_i(t)+\varepsilon]$.
On this interval the external field \eqref{netexternalfield}
is strictly convex (we use that the measures $\mu_j$ with $j \neq i$
are supported outside $[x_i(t)-\varepsilon,x_i(t) + \varepsilon]$).
The convexity implies that the support of $\mu_i$ is an interval, see e.g.\
\cite[Theorem IV 1.10]{SaffTotik}.

(b) The real analyticity of the density of $\mu_i$ in the interior of
the support follows from \cite{DKM},
since the effective external field \eqref{netexternalfield} is real analytic
on the support of $\mu_i$.

The convexity of \eqref{netexternalfield} implies that the density of $\mu_i$
does not vanish in the interior of the support, and has
square root behavior at the endpoints,  see
e.g.\ \cite[Lemma 3.5]{Claeys1}.
\end{proof}

\subsection{Varying $n$}

In the situation of Theorem \ref{theorem:limitingdistribution} we have for each
finite $n$, the number $n_{k,l}$ of paths going from $a_k$ to $b_l$. The finite
$n$ transition numbers are
\[ t_{k,l}^{(n)} = \frac{n_{k,l}}{n}  \]
and in the limit we have
\begin{equation} \label{limittkl}
    \lim_{n \to \infty} t_{k,l}^{(n)} = t_{k,l} \qquad
    \textrm{for } k=1, \ldots, p, \, l=1, \ldots, q.
    \end{equation}

The equilibrium problem depends on the transition numbers by means
of the normalizations \eqref{totalmassmui}. For a finite $n$,
we use the equilibrium problem for a vector of measures $(\mu_1, \ldots, \mu_{p+q-1})$ that
have total masses
\begin{equation} \label{n-normalization}
    \int d\mu_i = t_{k(i),l(i)}^{(n)}, \qquad \textrm{for } i =1, \ldots, p+q-1
    \end{equation}
instead of \eqref{totalmassmui}. Then the minimizer $(\mu_1^{(n)}, \ldots, \mu_{p+q-1}^{(n)})$
will also depend on $n$.

Because of Theorem \ref{theorem:energyminimizerexists}
and \eqref{limittkl}, there is a $T_{0} > 0$,
so that for every $T < T_{0}$, each $\mu_i^{(n)}$ is supported on
an interval $[\alpha_i^{(n)}, \beta_i^{(n)}]$ (depending on $n$)
so that parts (a) and (b) of Theorem \ref{theorem:energyminimizerexists} hold.
As $n \to \infty$, we have that $\mu_i^{(n)} \to \mu_i$  and
\[ \alpha_i^{(n)} \to \alpha_i, \qquad \beta_i^{(n)} \to \beta_i. \]
for every $i$.

In the steepest descent analysis that follows we will fix a large enough $n$.
We will work with the $n$-dependent measures $\mu_i^{(n)}$
and intervals $[\alpha_i^{(n)}, \beta_i^{(n)}]$, but for ease of notation
we will usually not write the superscript $(n)$. We trust that this
will not lead to confusion. However, we will write $t_{k,l}^{(n)}$.
A property that will be used a number of times is that
\begin{equation} \label{ntlkisinteger}
    n t_{k,l}^{(n)} = n_{k,l} \qquad \textrm{is an integer.}
    \end{equation}

\section{Riemann surface, $\xi$-functions, $\lam$-functions}
\label{section:xifunctions}

\subsection{Variational conditions}
\label{subsection:variational}

As said before, we take a large $n$ and consider the vector
equilibrium problem with normalization \eqref{n-normalization}
and we assume that $T$ is sufficiently small so that
Theorem \ref{theorem:energyminimizerexists} applies.
So the measure $\mu_i$ is supported on the interval $[\alpha_i, \beta_i]$.

The variational conditions associated with
the vector equilibrium problem \eqref{energyminimization:mixedangelesco}
are as follows. For each $i = 1, \ldots, p+q-1$, there is a constant
$L_i\in \mathbb R$ so that
\begin{equation}
\label{eulerlagrange:inequality} 2\sum_{j}
a_{i,j}\int \log\frac{1}{|x-y|}\, d\mu_j(y) +
    \frac{1}{T} V_i(x) \left\{\begin{array}{l} = L_i, \quad x\in [\alpha_i, \beta_i],\\
\geq L_i,\quad x\in\mathbb R \setminus [\alpha_i,\beta_i]. \end{array}\right.
\end{equation}

We use $F_i$ to denote the Cauchy transform of the measure $\mu_i$,
\begin{equation}\label{defcauchytransform}
    F_i(z) := \int_{\alpha_i}^{\beta_i} \frac{1}{z-x}\, d\mu_i(x),
\end{equation}
for $i=1,\ldots,p+q-1$. The function $F_i(z)$ is analytic on
$\mathbb C \setminus [\alpha_i,\beta_i]$.
By taking the derivative of \eqref{eulerlagrange:inequality} and using the
Sokhotski-Plemelj formula it follows that
\begin{equation}\label{eulerlagrange:equalityderivative}
    -F_{i,+}(x)-F_{i,-}(x) - 2\sum_{j\neq i} a_{i,j}F_{j}(x) + \frac{1}{T} V_i'(x) = 0,
    \quad x\in [\alpha_i, \beta_i].
\end{equation}

\begin{lemma} \label{lemma:eulerlagrange:strictinequality}
The variational inequality \eqref{eulerlagrange:inequality} is strict for $x\in
[\beta_{i+1},\alpha_{i})\cup (\beta_{i},\alpha_{i-1}]$, where
we put $\alpha_0 = +\infty$ and $\beta_{p+q} = -\infty$.
\end{lemma}

\begin{proof} On both gaps
$(\beta_{i+1},\alpha_{i})$ and $(\beta_{i},\alpha_{i-1})$, the left-hand side
of \eqref{eulerlagrange:inequality} is a real analytic function of $x$
whose first  derivative is
\begin{equation} \label{eulerlagrange:derivative}
-2\sum_{j} a_{i,j}F_{j}(x) + \frac{1}{T} V_i'(x)
\end{equation}
and whose second derivative is
\begin{equation}\label{eulerlagrange:secondderivative}
    2\sum_{j} a_{i,j}\int_{\alpha_i}^{\beta_i} \frac{1}{(x-y)^2}\, d\mu_j(y) +
    \frac{1}{T} V_i''(x).
\end{equation}
Each term in \eqref{eulerlagrange:secondderivative} is positive
and so the left-hand side of \eqref{eulerlagrange:inequality}
is strictly convex on both $(\beta_{i+1},\alpha_{i})$ and $(\beta_{i},\alpha_{i-1})$,
which proves the lemma.
\end{proof}

\subsection{Riemann surface}
\label{subsection:xifunctions}

We construct a Riemann surface $\mathcal R$ as follows, compare with \cite{Duits}.
The Riemann surface has $p+q$ sheets which we denote by $\mathcal R_j$, $j=1, \ldots, p+q$.
Each sheet is associated with a vertex of the graph, that is, with either
a starting point $a_k$ or an ending point $b_l$. We choose the numbering
so that $\mathcal R_k$ is associated with $a_k$ for $k=1, \ldots, p$,
and $\mathcal R_{p+l}$ is associated with $b_l$ for $l=1, \ldots, q$.

\begin{figure}[t]
\begin{center}
   \setlength{\unitlength}{1truemm}
   \begin{picture}(100,70)(-5,2)

       \put(10,45){\line(1,0){55}}
       \put(20,55){\line(1,0){55}}
       \put(10,45){\line(1,1){10}}
       \put(65,45){\line(1,1){10}}
       \thicklines
       \put(55,50){\line(1,0){8}}
       \thinlines
       \put(55,50){\thicklines\circle*{1}}
       \put(63,50){\thicklines\circle*{1}}
       \put(52,52){$\alpha_1$}
       \put(63,52){$\beta_1$}
       \put(75,50){$\xi_{1}$}

       \put(10,30){\line(1,0){55}}
       \put(20,40){\line(1,0){55}}
       \put(10,30){\line(1,1){10}}
       \put(65,30){\line(1,1){10}}
       \thicklines
       \put(20,35){\line(1,0){8}}
       \put(35,35){\line(1,0){10}}
       \thinlines
       \put(20,35){\thicklines\circle*{1}}
       \put(28,35){\thicklines\circle*{1}}
       \put(35,35){\thicklines\circle*{1}}
       \put(45,35){\thicklines\circle*{1}}
       \put(17,37){$\alpha_3$}
       \put(28,37){$\beta_3$}
       \put(32,37){$\alpha_2$}
       \put(45,37){$\beta_2$}
       \put(75,35){$\xi_{2}$}

       \put(10,15){\line(1,0){55}}
       \put(20,25){\line(1,0){55}}
       \put(10,15){\line(1,1){10}}
       \put(65,15){\line(1,1){10}}
       \thicklines
       \put(35,20){\line(1,0){10}}
       \put(55,20){\line(1,0){8}}
       \thinlines
       \put(35,20){\thicklines\circle*{1}}
       \put(45,20){\thicklines\circle*{1}}
       \put(55,20){\thicklines\circle*{1}}
       \put(63,20){\thicklines\circle*{1}}
       \put(32,22){$\alpha_2$}
       \put(45,22){$\beta_2$}
       \put(52,22){$\alpha_1$}
       \put(63,22){$\beta_1$}
       \put(75,20){$\xi_{3}$}

       \put(10,0){\line(1,0){55}}
       \put(20,10){\line(1,0){55}}
       \put(10,0){\line(1,1){10}}
       \put(65,0){\line(1,1){10}}
       \thicklines
       \put(20,5){\line(1,0){8}}
       \thinlines
       \put(20,5){\thicklines\circle*{1}}
       \put(28,5){\thicklines\circle*{1}}
       \put(17,7){$\alpha_3$}
       \put(28,7){$\beta_3$}
       \put(75,5){$\xi_{4}$}

   \end{picture}
   \caption{Four-sheeted Riemann surface for the example in Figure \ref{fig:graph1}.}
   \label{fig:Riemannsurface4x4}
\end{center}
\end{figure}
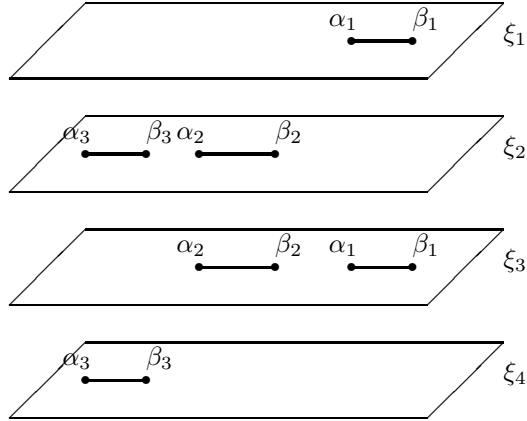

Recall that we use $i$ as a label for the edges of the graph, and we write
$k =k(i)$ and $l=l(i)$ if $i$ labels the edge $(a_k,b_l)$.
Then the $p+q$ sheets are defined as
\begin{align}
    \mathcal R_k & = \overline{\mathbb C}\setminus\bigcup_{i :\, k(i)=k} [\alpha_i, \beta_i],
    \quad k=1,\ldots,p,\\
    \mathcal R_{p+l} & = \overline{\mathbb C}\setminus\bigcup_{i: \, l(i)=l} [\alpha_i, \beta_i] ,\quad
    l=1,\ldots,q.
\end{align}
The sheets are connected as follows. For each $i=1, \ldots, p+q-1$ we have
that $\mathcal R_{k(i)}$ is  connected to sheet $\mathcal R_{p+l(i)}$ along the interval $[\alpha_i,\beta_i]$
in the usual crosswise manner.
See Figure
\ref{fig:Riemannsurface4x4} for a picture of the Riemann surface for the
example in Figure \ref{fig:graph1}. Note that the cuts $[\alpha_i,\beta_i]$
are always between a sheet $\mathcal R_k$ which is among the first $p$ sheets
and a sheet $\mathcal R_{p+l}$ which is among the last $q$ sheets.
The cuts strictly move to the left if we go from one sheet to the next among the
first $p$ sheets, and similarly, among the last $q$ sheets.

The Riemann surface depends on $n$, since the endpoints $\alpha_i^{(n)}$ and $\beta_i^{(n)}$
depend on $n$. The $\xi$- and $\lambda$-functions that we define from the Riemann surface
will also depend on $n$. We will not indicate the $n$-dependence, as already mentioned
before.

The Riemann surface $\mathcal R$ has $p+q$ sheets and $2(p+q-1)$ simple branch
points. Therefore by Hurwitz's formula (see e.g. \cite{Mir}) its genus
$g$ satisfies
 \[ 2g-2 = -2(p+q) + 2(p+q-1) = -2 \]
 so that $g=0$.
The fact that the genus is zero will be helpful in the construction of the
global parametrix in Section~\ref{subsection:globalparametrix}.

\subsection{$\xi$-functions}

We define the $\xi$-functions as follows. Recall that $F_i$ is
given by \eqref{defcauchytransform}.
\begin{align}\label{xifunctions:mixedangelesco:a}
    \xi_{k}(z) & = -\sum_{i :\, k(i)=k} F_i(z)+\frac{1}{Tt} (z-a_k),& k=1,\ldots,p,\\
    \label{xifunctions:mixedangelesco:b}
    \xi_{p+l}(z) & = \sum_{i: \, l(i)=l} F_i(z)-\frac{1}{T(1-t)} (z-b_l), & l=1,\ldots,q.
\end{align}

We consider $\xi_j$ as an analytic function on the $j$th sheet $\mathcal R_j$
of the Riemann surface with a pole at infinity. Moreover, these functions
define a global meromorphic function on $\mathcal R$ as the following result
shows.

\begin{theorem} \label{theorem:xifunctionscompatible}
Consider the function $\xi_{j}(z)$ on the $j$th sheet $\mathcal R_j$,
$j=1,\ldots,p+q$. Then these functions are compatible along the cuts $[\alpha_i,\beta_i]$
of the Riemann surface in the sense that
\begin{equation} \label{eulerlagrange:equalityxifunctions}
\left\{\begin{array}{l} \xi_{k(i),+}(x) = \xi_{p+l(i),-}(x),\\
\xi_{k(i),-}(x)=\xi_{p+l(i),+}(x),\end{array}\right. \quad x\in [\alpha_i, \beta_i],
\end{equation} for every $i=1,\ldots,p+q-1$. Hence the $\xi_j$-functions can
be extended to a global meromorphic function $\xi$ defined on the Riemann
surface $\mathcal R$.
\end{theorem}

\begin{proof} Fix $i\in\{1,\ldots,p+q-1\}$. On the interval $[\alpha_i,\beta_i]$ we have
by definition
\eqref{xifunctions:mixedangelesco:a}--\eqref{xifunctions:mixedangelesco:b} that
    \begin{multline} \label{xifunctionsjump}
        \xi_{k(i),+}(x)-\xi_{p+l(i),-}(x) \\
    = - \sum_{j :\, k(j)=k(i)} F_{j,+}(x) - \sum_{j: \, l(j)=l(i)} F_{j,-}(x)
    +\frac{1}{Tt(1-t)}(x-(1-t)a_{k(i)}-tb_{l(i)})) \\
    =  -F_{i,+}(x)-F_{i,-}(x)-
        \left(\sum_{j\neq i :\, k(j) = k(i) \textrm{ or } l(j)=l(i)} F_{j}(x)\right)
            +  \frac{1}{T} V'_i(x),
\end{multline}
where we used the definition \eqref{defxit}--\eqref{defexternalfield} of $V_i$.

If $j \neq i$ is such that $k(j) = k(i)$ or $l(j) = l(i)$ then $a_{i,j} = 1/2$.
For all other $j \neq i$ we have $a_{i,j} = 0$. Therefore
\[ \left(\sum_{j\neq i :\, k(j) = k(i) \textrm{ or } l(j)=l(i)} F_{j}(x)\right)
    = 2 \sum_{j \neq i} a_{i,j} F_j(x). \]
Using this in \eqref{xifunctionsjump} and recalling the variational equality
\eqref{eulerlagrange:equalityderivative}, we see that that
$\xi_{k(i),+}(x)=\xi_{p+l(i),-}(x)$ for  $x \in [\alpha_i,\beta_i]$.

The other
equality $\xi_{k(i),-}(x)=\xi_{p+l(i),+}(x)$ follows in exactly the same way.
\end{proof}

In the following two lemmas we collect some more properties
of the $\xi$-functions that will be needed in what follows.

\begin{lemma} \label{lemma:xiasymptotics}
The $\xi$-functions have the following behavior as $z \to \infty$
\begin{align}
\label{asymptotics:xia}
    \xi_{k}(z) & = \frac{1}{Tt} (z-a_k) - \frac{n_k}{n z} + O\left(\frac{1}{z^2}\right),
        \qquad k=1, \ldots, p, \\
\label{asymptotics:xib} \xi_{p+l}(z) & =
    - \frac{1}{T(1-t)} ( z-b_l ) + \frac{m_l}{n z} +O\left(\frac{1}{z^2}\right),
        \qquad l=1, \ldots, q.
\end{align}
\end{lemma}
\begin{proof}
From \eqref{defcauchytransform} and \eqref{n-normalization} it follows that
\begin{equation} \label{Fiatinfinity}
    F_i(z) = \frac{\int d\mu_i}{z}  +  O\left(\frac{1}{z^2}\right)
    = \frac{t_{k(i),l(i)}^{(n)}}{z} + O\left(\frac{1}{z^2}\right)
    \end{equation}
    as $z \to \infty$. Recall that we work with $n$-dependent transition numbers.
Since
\[ \sum_{i: \, k(i) = k} t_{k(i),l(i)}^{(n)} = \frac{n_k}{n},
    \qquad \sum_{i:\, l(i) = l} t_{k(i),l(i)}^{(n)} = \frac{m_l}{n} \]
the asymptotic behaviors \eqref{asymptotics:xia} and \eqref{asymptotics:xib}
follow immediately from the definitions
\eqref{xifunctions:mixedangelesco:a}--\eqref{xifunctions:mixedangelesco:b}.
\end{proof}

\begin{lemma}\label{lemma:integralalongcut1}
For any $i=1,\ldots,p+q-1$ we have
\begin{align}
\int_{\ccal_i} \xi_{k(i)}(z)\, dz &= -2\pi {\rm i} \, t_{k(i),l(i)}^{(n)},\\
\int_{\ccal_i} \xi_{p+l(i)}(z)\, dz &= 2\pi{\rm i} \, t_{k(i),l(i)}^{(n)},\\
\int_{\ccal_i} \xi_{j}(z)\, dz &= 0,\qquad\qquad j\not\in\{k(i),p+l(i)\},
\end{align}
where $\ccal_i$ denotes a counterclockwise contour surrounding the interval $[\alpha_i,\beta_i]$
and not enclosing any point of the other intervals $[\alpha_j,\beta_j]$, $j\neq i$.
\end{lemma}

\begin{proof}
From \eqref{defcauchytransform}, \eqref{Fiatinfinity} and the
definition of $\ccal_i$ it follows that
\[ \int_{\ccal_i} F_i(z) \, dz = 2 \pi {\rm i} \, t_{k(i),l(i)}^{(n)} \]
and
\[ \int_{\ccal_i} F_j(z) \, dz = 0 \qquad \textrm{if } j \neq i. \]
The lemma then follows from the definitions
\eqref{xifunctions:mixedangelesco:a}--\eqref{xifunctions:mixedangelesco:b}.
\end{proof}

\begin{remark}
Our definition of the $\xi$-functions differs slightly from the one
used in \cite{BK2}. If $\tilde{\xi}_j$ denote the
$\xi$-functions in \cite{BK2} then we have
\[ \xi_j(z)= \tilde{\xi}_j(z) -\frac{z}{T(1-t)}, \qquad j=1,\ldots,p+q. \]
In the present form the formulae are more symmetric.
\end{remark}

\subsection{$\lambda$-functions}
\label{subsection:lambdafunctions}

We define the $\lambda_j$-functions as
\begin{align}\label{lamfunctions:mixedangelesco:a}
\lambda_{k}(z) & = c_{k} +\int_{\beta_{i_{k}}}^{z}\xi_{k}(s)\, ds, & k=1,\ldots,p,\\
\label{lamfunctions:mixedangelesco:b}
\lambda_{p+l}(z) & = c_{p+l}+\int_{\beta_{\tilde{i}_{l}}}^{z}\xi_{p+l}(s)\, ds, & l=1,\ldots,q,
\end{align}
where
\[ i_k := \min\{i \mid k(i)=k\}, \quad \text{and} \quad
    \tilde{i}_l := \min\{i \mid l(i)=l\}, \]
and where the path of integration in the integrals
\eqref{lamfunctions:mixedangelesco:a} and \eqref{lamfunctions:mixedangelesco:b}
lies in $\mathbb C \setminus (-\infty,\beta_{i_k})$ and
$\mathbb C \setminus (-\infty,\beta_{\tilde{i}_l})$, respectively.
The functions $\lambda_k(z)$ and $\lambda_{p+l}(z)$ are defined with a branch cut along
the intervals $(-\infty,\beta_{i_k}]$ and $(-\infty,\beta_{\tilde{i}_l}]$,
respectively.


We choose the constants $c_{j}$  in
in the following way.

\begin{lemma}
We can (and do) choose the constants $c_j$ in
\eqref{lamfunctions:mixedangelesco:a}--\eqref{lamfunctions:mixedangelesco:b}
in such a way that
\begin{equation}\label{realpartlambdas1}
    \Re(\lambda_{k(i),+}(\beta_i)) = \Re(\lambda_{p+l(i),+}(\beta_i)),
\end{equation}
for every $i=1,\ldots,p+q-1$.
\end{lemma}

\begin{proof}
We use the fact that the graph
$G=(V,E,t)$ is a tree. We can iteratively \lq undress\rq\ this tree as follows.
We start with $G_1=G$. Next we choose a leaf vertex $v$ and set
$$G_2=G_1\setminus v,$$ i.e., $G_2$ is the tree obtained by removing the leaf
and its corresponding edge from the tree $G_1$. We iteratively repeat this
procedure and obtain in this way a chain of nested trees
\begin{equation}\label{chainofnestedtrees}
G = G_1\supset G_2\supset \cdots \supset G_{|V|},
\end{equation}
where each $G_{i}$ is obtained from $G_{i-1}$ by removing one leaf. Obviously the
last non-empty tree $G_{|V|}$ in this chain consists of a single vertex $v_j$.

We freely choose the corresponding constant $c_{j}$. Next we use induction
on $k=|V|-1,\ldots,2,1$, to fill in all the remaining constants $c_{j}$ so that
each time \eqref{realpartlambdas1} is satisfied. This is possible since
$G_{i-1}\setminus G_{i}$ consists of a single vertex $v_j$, which is a leaf of
$G_{i-1}$, and hence we have exactly one condition \eqref{realpartlambdas1} to
fix the integration constant $c_j$ of this leaf. Thus we see that the
conditions \eqref{realpartlambdas1} can indeed be imposed on the constants
$c_{j}$ in
\eqref{lamfunctions:mixedangelesco:a}--\eqref{lamfunctions:mixedangelesco:b}.
\end{proof}

Properties of the $\lambda$-functions that we will need are stated
in the following lemmas. The first is a reformulation
of the variational conditions
\eqref{eulerlagrange:inequality}.

\begin{lemma} \label{lem:eulerlagrangerestatement}
We have for $i=1,\ldots,p+q-1$,
\begin{equation} \label{eulerlagrange:lambda2}
\Re(\lambda_{k(i),\pm}(x)-\lambda_{p+l(i),\pm}(x))
\left\{\begin{array}{l} = 0,\quad x\in [\alpha_i,\beta_i],\\
\geq 0,\quad x\in\mathbb R \setminus [\alpha_i,\beta_i]. \end{array}\right.
\end{equation}
Strict inequality in \eqref{eulerlagrange:lambda2} holds for
$x \in [\beta_{i+1}, \alpha_i) \cup (\beta_i, \alpha_{i-1}]$ (where $\alpha_0=+\infty$ and $\beta_{p+q}=-\infty$).
\end{lemma}
\begin{proof} Observe that by
\eqref{xifunctions:mixedangelesco:a}--\eqref{xifunctions:mixedangelesco:b} and
\eqref{lamfunctions:mixedangelesco:a}--\eqref{lamfunctions:mixedangelesco:b} we
have that
\begin{align}
\label{lambda:explicit1}
    \lambda_{k}(z) &=  \sum_{i :\, k(i)=k}  \int \log\frac{1}{z-x}\, d\mu_i(x)
    +\frac{1}{2Tt} (z-a_k)^2 + \tilde{c}_k,\\
\label{lambda:explicit2}
    \lambda_{p+l}(z) & = -\sum_{i : \, l(i)=l} \int \log\frac{1}{z-x}\, d\mu_i(x)
    -\frac{1}{2T(1-t)} (z-b_l)^2 + \tilde{c}_{p+l},
\end{align}
for certain real constants $\tilde{c}_j$, $j=1,\ldots,p+q$.
Therefore, for $i=1,\ldots,p+q-1$,
\begin{multline}\label{eulerlagrange:lambda0}
    \Re(\lambda_{k(i)}(z) - \lambda_{p+l(i)}(z)) \\
    = 2\int \log\frac{1}{|z-x|}\, d\mu_i(x)
    + \sum_{j\neq i :\, k(j)=k(i) \textrm{ or } l(j)=l(i)} \int \log\frac{1}{|z-x|}\, d\mu_j(x) \\
    +\frac{1}{T} \Re V_i(z) + \frac{1}{2T}(a_{k(i)} - b_{l(i)})^2 +\tilde{c}_{k(i)}-\tilde{c}_{p+l(i)},
\end{multline}
where we used the definition \eqref{defexternalfield} of $V_i$.
Then by the variational conditions \eqref{eulerlagrange:inequality}
we have
\begin{equation}\label{eulerlagrange:lambda1}
    \Re(\lambda_{k(i),\pm}(x) - \lambda_{p+l(i),\pm}(x))
    \geq L_i + \frac{1}{2T}(a_{k(i)} - b_{l(i)})^2 +\tilde{c}_{k(i)}-\tilde{c}_{p+l(i)}
\end{equation}
for $x \in \mathbb R$ with equality  for $x \in [\alpha_i,\beta_i]$.
The  constant
in the right-hand side of \eqref{eulerlagrange:lambda1} is equal to zero
because the constants $c_j$ are chosen so that \eqref{realpartlambdas1} holds.

The strict inequality for
$x \in [\beta_{i+1}, \alpha_i) \cup (\beta_i, \alpha_{i-1}]$
is a consequence of Lemma \ref{lemma:eulerlagrange:strictinequality}.
\end{proof}

\begin{lemma}
As $z\to\infty$ we have that
\begin{align}
\label{asymptotics:lambdaa} \lambda_{k}(z) & =
    \frac{1}{2Tt} (z-a_k)^2 - \frac{n_k}{n}\log z+ \tilde{c}_{k}+O\left(\frac{1}{z}\right),
        \quad k =1, \ldots, p,  \\
\label{asymptotics:lambdab} \lambda_{p+l}(z) & =
     -\frac{1}{2T(1-t)}(z-b_l)^2 +\frac{m_l}{n}\log z+\tilde{c}_{p+l}+O\left(\frac{1}{z}\right),
    \quad l=1, \ldots, q.
\end{align}
\end{lemma}
\begin{proof}
This follows from the definitions \eqref{lamfunctions:mixedangelesco:a}
and \eqref{lamfunctions:mixedangelesco:b}
and the asymptotic behavior \eqref{asymptotics:xia} and \eqref{asymptotics:xib}
of the $\xi$-functions.
\end{proof}

The next lemma will be a consequence of Lemma \ref{lemma:integralalongcut1}.

\begin{lemma}\label{lemma:integralalongcut2}
For  $k=1,\ldots,p$, we have that
\begin{equation}\label{jumps:negativerealline1}
    \exp(n(\lambda_{k,+}(x)-\lambda_{k,-}(x))) = 1,
        \qquad \textrm{for } x\in\mathbb R \setminus \bigcup_{i:\, k(i) =k} [\alpha_i,\beta_i].
\end{equation}

For $l=1, \ldots, q$, we have that
\begin{equation}\label{jumps:negativerealline2}
    \exp(n(\lambda_{p+l,+}(x)-\lambda_{p+l,-}(x))) = 1, \qquad
        \textrm{for } x \in \mathbb R \setminus \bigcup_{i:\, l(i) =l} [\alpha_i,\beta_i].
\end{equation}
\end{lemma}

\begin{proof} Fix $k=1, \ldots, p$ and let
$x\in\mathbb R \setminus \bigcup_{i:\, k(i) =k} [\alpha_i,\beta_i]$.
By definition of $\lambda_k$ and by contour deformation we have that
\[ \lambda_{k,+}(x)-\lambda_{k,-}(x) = \int_{\ccal} \xi_{k}(z)\, dz
\]
where $\ccal$ is some closed contour surrounding some of the intervals
$[\alpha_i,\beta_i]$. From Lemma \ref{lemma:integralalongcut1} it follows that
each of the enclosed intervals gives a contribution to the integral of the form
$\pm 2\pi{\rm i} \, t_{k,l}^{(n)}$. Since each $t_{k,l}^{(n)}$ is a rational number with
denominator $n$, indeed $t_{k,l}^{(n)} = n_{k,l}/n$, we conclude that
$n(\lambda_{k,+}(x)-\lambda_{k,-}(x))$ is a multiple of $2\pi{\rm i}$, and so
we obtain \eqref{jumps:negativerealline1}.

The proof of \eqref{jumps:negativerealline2} is similar.
\end{proof}

Lemma \ref{lemma:integralalongcut2} will be important for the steepest descent
analysis in Section~\ref{section:normalizationatinfinity}.

\section{First transformation of the RH problem: Normalization at infinity}
\label{section:normalizationatinfinity}

In the following sections we describe the steepest descent analysis $Y\mapsto
X\mapsto T\mapsto S\mapsto R$ of the RH problem \ref{RHP:original}. Throughout
the steepest descent analysis the following simple lemma will be repeatedly
used.

\begin{lemma}\label{lemma:transformation:jumps}
Assume that the matrix function $Y(z)$ satisfies the jump condition $Y_{+}(x) =
Y_-(x)J(x)$ for $x\in\mathbb R$. Let $A(z)$ and $B(z)$ be matrix functions with
$A(z)$ entire and $B(z)$ analytic in $\mathbb C\setminus\mathbb R$. Then
$$X(z) := A(z)Y(z)B(z)$$ satisfies the jump condition $$X_+(x) = X_-(x)\left(
B^{-1}_-(x)J(x)B_+(x) \right).$$
\end{lemma}

The point will be to choose appropriate transformation matrices $A(z)$ and
$B(z)$ in order to bring the RH problem \ref{RHP:original} to a simple
form.

Let $Y$ be the solution to the original RH problem \ref{RHP:original}. The
first transformation $Y\mapsto X$ serves to normalize the RH problem at
infinity. To this end we define $X=X(z)$ as
\begin{equation}\label{defX}
X(z) = L^{-n} Y(z)D(z)^n,
\end{equation}
where we define the diagonal matrices
\begin{multline}\label{defD}
    D(z) = \diag\left(\left(\exp\left(\lambda_{k}(z)-
        \frac{1}{2Tt}(z - a_k)^2 \right)\right)_{k=1}^p,\right.\\
    \left.\left(\exp\left(\lambda_{p+l}(z)+ \frac{1}{2T(1-t)} (z-b_l)^2\right)\right)_{l=1}^q\right),
\end{multline}
and
\begin{equation}\label{defL}
    L = \diag(\exp(\tilde{c}_{1}),\ldots,\exp(\tilde{c}_{p+q})),
\end{equation}
where the constants $\tilde{c}_j$ are as in
\eqref{asymptotics:lambdaa}--\eqref{asymptotics:lambdab}. From Lemma
\ref{lemma:transformation:jumps} it follows by straightforward calculations
that $X=X(z)$ satisfies the following RH problem.

\begin{rhp} \textrm{ }
\begin{itemize}
\item[\rm (1)] $X$ is analytic
in $\mathbb C \setminus \mathbb R$;
 \item[\rm (2)] For $x\in \mathbb R$, we have that
\begin{equation}\label{jumpsX1}
    X_{+}(x) = X_{-}(x) \begin{pmatrix}
    J_{1,1}(x) & J_{1,2}(x)\\
    0 & J_{2,2}(x)
\end{pmatrix}
\end{equation}
where the blocks $J_{1,1}$, $J_{1,2}$ and $J_{2,2}$
(of sizes $p \times p$, $p \times q$ and $q \times q$, respectively)
have the following form.
\begin{enumerate}
\item[\rm (a)] $J_{1,2}$ is a full $p \times q$ matrix with entries
\begin{equation} \label{jumpsX4}
   \left( J_{1,2}(x) \right)_{k,l} = \exp\left(n(\lambda_{p+l,+}(x)-\lambda_{k,-}(x))\right),
    \qquad x \in \mathbb R.
\end{equation}
\item[\rm (b)] Outside of the intervals $[\alpha_i,\beta_i]$,
$J_{1,1}$ and $J_{2,2}$ are identity matrices
\begin{equation}\label{jumpsX5}
J_{1,1}(x) = I_p, \qquad J_{2,2}(x) = I_q,
    \qquad x \in \mathbb R\setminus\bigcup_{i=1}^{p+q-1} [\alpha_i,\beta_i].
\end{equation}
\item[\rm (c)] On the interval $[\alpha_i,\beta_i]$,
$J_{1,1}$ and $J_{2,2}$ are diagonal matrices with  ones on the diagonal,
except for
\begin{equation} \label{jumpsX2}
    \left(J_{1,1}(x) \right)_{k(i),k(i)} = \exp\left(n(\lambda_{k(i),+}(x)-\lambda_{k(i),-}(x))\right),
    \qquad x \in (\alpha_i, \beta_i), \\
\end{equation}
and
\begin{equation} \label{jumpsX3}
    \left(J_{2,2}(x)\right)_{l(i),l(i)} = \exp\left(n(\lambda_{p+l(i),+}(x)-\lambda_{p+l(i),-}(x))\right),
    \qquad x \in (\alpha_i,\beta_i).
\end{equation}
\end{enumerate}
\item[\rm (3)] As $z\to\infty$, we have that
\begin{equation}\label{asymptoticconditionX} X(z) =
    I_{p+q}+O(1/z).
\end{equation}
\end{itemize}
\end{rhp}
Here we used Lemma \ref{lemma:integralalongcut2} to see that the diagonal
entries of \eqref{jumpsX5}, as well as most of the diagonal entries of
\eqref{jumpsX2}--\eqref{jumpsX3} are equal to $1$. On the other hand, we used
the asymptotic conditions
\eqref{asymptotics:lambdaa}--\eqref{asymptotics:lambdab} and
\eqref{defX}--\eqref{defL} to see that the RH problem for $X$ is normalized at
infinity in the sense of \eqref{asymptoticconditionX}.

Let us illustrate the jump matrices for the example with $p=q=2$ as in Figures
\ref{fig:graph1} and \ref{fig:Riemannsurface4x4}. In that case the jump
conditions are written as
\begin{equation*}
X_{+} = X_{-}
\begin{pmatrix}
\exp(n(\lambda_{1,+}-\lambda_{1,-})) & 0 & \exp(n(\lambda_{3,+}-\lambda_{1,-})) &\exp(n(\lambda_{4}-\lambda_{1,-})) \\
0 & 1 & \exp(n(\lambda_{3,+}-\lambda_{2\phantom{,-}})) & \exp(n(\lambda_{4}-\lambda_{2\phantom{,-}})) \\
0 & 0 & \exp(n(\lambda_{3,+}-\lambda_{3,-})) &0 \\
0 & 0 & 0 & 1
\end{pmatrix}
\end{equation*}
on the interval $[\alpha_1,\beta_1]$,
\begin{equation*}
X_{+} = X_{-}
\begin{pmatrix}
1 & 0 &  \exp(n(\lambda_{3,+}-\lambda_{1\phantom{,-}})) &\exp(n(\lambda_{4}-\lambda_{1\phantom{,-}})) \\
0 & \exp(n(\lambda_{2,+}-\lambda_{2,-})) & \exp(n(\lambda_{3,+}-\lambda_{2,-})) & \exp(n(\lambda_{4}-\lambda_{2,-}))\\
0 & 0 & \exp(n(\lambda_{3,+}-\lambda_{3,-})) & 0 \\
0 & 0 & 0 & 1
\end{pmatrix}
\end{equation*}
on the interval $[\alpha_2,\beta_2]$,
\begin{equation*}
X_{+} = X_{-}
\begin{pmatrix}
1 & 0 & \exp(n(\lambda_{3}-\lambda_{1\phantom{,-}})) &\exp(n(\lambda_{4,+}-\lambda_{1\phantom{,-}})) \\
0 & \exp(n(\lambda_{2,+}-\lambda_{2,-})) & \exp(n(\lambda_{3}-\lambda_{2,-})) & \exp(n(\lambda_{4,+}-\lambda_{2,-})) \\
0 & 0 & 1 & 0 \\
0 & 0 & 0 & \exp(n(\lambda_{4,+}-\lambda_{4,-})) \\
\end{pmatrix}
\end{equation*}
on the interval $[\alpha_3,\beta_3]$, and
\begin{equation*}
X_{+} = X_{-}
\begin{pmatrix}
1 & 0 & \exp(n(\lambda_{3}-\lambda_{1})) &\exp(n(\lambda_{4}-\lambda_{1})) \\
0 & 1 & \exp(n(\lambda_{3}-\lambda_{2})) & \exp(n(\lambda_{4}-\lambda_{2})) \\
0 & 0 & 1 & 0 \\
0 & 0 & 0 & 1 \\
\end{pmatrix}
\end{equation*}
on $\mathbb R\setminus \bigcup_{i=1}^3 [\alpha_i,\beta_i]$.

\section{Second transformation of the RH problem: Gaussian elimination and opening of global lenses}
\label{section:globallenses}

On the interval $[\alpha_i, \beta_i]$ we have that $\lambda_{p+l(i),+} - \lambda_{k(i),-}$
is a purely imaginary constant. So the entry
\begin{equation} \label{J12constant}
    \left(J_{1,2}(x) \right)_{k(i),l(i)}, \qquad x \in [\alpha_i, \beta_i]
    \end{equation}
is constant with modulus one.
It would be an ideal situation if, except for \eqref{J12constant},
all entries of the matrix $J_{1,2}(x)$ for $x \in \mathbb R$, are exponentially decaying
as $n \to \infty$. That would happen if for every $k=1, \ldots, p$, and $l = 1, \ldots, q$,
we have
\begin{enumerate}
\item[(a)] if $t_{k,l}^{(n)} = 0$ then
\[ \Re \lambda_{p+l,+}(x) < \Re \lambda_{k,-}(x), \qquad x \in \mathbb R, \]
\item[(b)] if $t_{k,l}^{(n)} > 0$ then
\[ \Re \lambda_{p+l,+}(x) < \Re \lambda_{k,-}(x), \qquad x \in \mathbb R \setminus [\alpha_i,\beta_i], \]
where $i = k+l-1$.
\end{enumerate}

However, this will not happen in general, and that is why we need
a second transformation $X \mapsto T$ in the steepest descent analysis
of the RH problem in which a number of unwanted entries of the jump matrices are eliminated.
In particular those entries of $J_{1,2}(x)$ that could potentially be exponentially
increasing as $n \to \infty$.

To this end we will open
global lenses \cite{ABK} and apply Gaussian elimination inside each of the
lenses. This construction will be systematic and may have interest in its
own right. The proof that appropriate lenses exist will be a consequence of the
maximum principle for subharmonic functions.

The opening of global lenses can be conveniently described in terms of the
right-down path in Proposition~\ref{lemma:latticepath}. We start at
the top left entry $(1,1)$ of the matrix of transition numbers $(t_{k,l}^{(n)})$ and
walk along the path until we
arrive at the bottom right entry $(p,q)$. During this walk, we will open global
lenses in an appropriate way. The precise action to perform depends on whether
we are taking a vertical (down) or a horizontal (right) step along the path.

\subsection{Construction of global lenses: Vertical step}
\label{subsection:verticalstep}

First we construct the global lenses for a vertical step along the lattice
path. Thus assume that $i\in\{1,\ldots,p+q-2\}$ is such that
\begin{align}
\label{verticalstep:condition1}
(k(i),l(i)) &= (k,l),\\
\label{verticalstep:condition2} (k(i+1),l(i+1)) &= (k+1,l),
\end{align}
for certain $k=1,\ldots,p-1$ and $l=1,\ldots,q$.

From \eqref{asymptotics:lambdaa} we obtain the asymptotic behavior
\begin{align}
 \lambda_{k+1}(z)-\lambda_k(z)
\label{asymptotics:differencelambdas}= & \frac{1}{Tt} (a_k-a_{k+1})z + O(\log |z|),
\end{align}
as $z \to \infty$.

We also note that $\Re \lambda_k$ and $\Re \lambda_{k+1}$ are well-defined
and continuous on $\mathbb C$. Indeed, we have by \eqref{lambda:explicit1}
that
\begin{align} \label{Realpartoflambdak}
 \Re(\lambda_{k}(z)) & = \sum_{j:\,  k(j)=k} \int \log\frac{1}{|z-x|} d\mu_j(x)
    + \frac{1}{2Tt} \Re (z-a_k)^2 + \tilde{c}_k, \\
    \label{Realpartoflambdak+1}
    \Re(\lambda_{k+1}(z)) & = \sum_{j:\,  k(j)=k+1}     \int \log\frac{1}{|z-x|}d\mu_j(x)   +
    \frac{1}{2Tt} \Re (z-a_{k+1})^2 + \tilde{c}_{k+1},
\end{align}
for certain constants $\tilde{c}_k$ and $\tilde{c}_{k+1}$.

The representations \eqref{Realpartoflambdak} and \eqref{Realpartoflambdak+1}
also show the following.

\begin{lemma} \label{lemma:subsuperharmonisch}
The function $z \mapsto \Re(\lambda_{k+1}(z)-\lambda_k(z))$ is superharmonic
on $\mathbb C \setminus \bigcup_{j:\,  k(j)=k} [\alpha_j,\beta_j]$, and subharmonic on
$\mathbb C\setminus\bigcup_{j:\,  k(j)=k+1} [\alpha_j,\beta_j]$.
\end{lemma}

\begin{proof} It is a standard fact from potential theory
that any function of the form $z\mapsto\int \log\frac{1}{|z-x|}\, d\mu(x)$,
with $\mu$ a positive measure with compact support, is superharmonic on
$\mathbb C$ and harmonic on $\mathbb C \setminus \supp(\mu)$, see e.g.\
\cite[Chapter 0]{SaffTotik}. Thus by \eqref{Realpartoflambdak} and \eqref{Realpartoflambdak+1}
we have that $z \mapsto \Re \lambda_{k+1}(z)$ is
superharmonic on $\mathbb C$ and harmonic on $\mathbb C \setminus \bigcup_{j:\,  k(j)=k+1} [\alpha_j,\beta_j]$,
while $z \mapsto - \Re \lambda_k(z)$ is subharmonic on $\mathbb C$
and harmonic on $\mathbb C \setminus \bigcup_{j:\,  k(j)=k} [\alpha_j,\beta_j]$

Since the two sets $ \bigcup_{j:\,  k(j)=k} [\alpha_j,\beta_j]$
and $\bigcup_{j:\,  k(j)=k+1} [\alpha_j,\beta_j]$ are disjoint,
the lemma follows.
 \end{proof}

We next define the open sets $\Omega_+,\Omega_-\subset\mathbb C$ as follows
\begin{align}\label{def:omegaplus}
    \Omega_+ & := \{z\in\mathbb C \mid  \Re(\lambda_{k+1}(z)-\lambda_k(z)) > 0\} \\
    \label{def:omegamin}
    \Omega_- &:= \{z\in\mathbb C \mid  \Re(\lambda_{k+1}(z)-\lambda_k(z))<0\}.
\end{align}
We also denote
\begin{align}\label{def:omegaplusinfty}
    \Omega_{+,\infty} & := \{z\in\Omega_+\mid  \exists\textrm{ a connected path in $\Omega_+$ from $z$ to $\infty$}\}\\
    \label{def:omegamininfty}
    \Omega_{-,\infty} & := \{z\in\Omega_-\mid  \exists\textrm{ a connected path in $\Omega_-$ from $z$ to $\infty$}\}.
\end{align}
In other words, $\Omega_{+,\infty}$ is the union of the unbounded connected
components of $\Omega_{+}$, and similarly for $\Omega_{-,\infty}$.

The open sets $\Omega_{+,\infty}$, $\Omega_{-,\infty}$ satisfy the following
properties.

\begin{lemma}\label{lemma:infinitecomponents}
\begin{enumerate}
\item[\rm (a)] For each $\varepsilon > 0$
there exists $R > 0$ so that
\begin{align} \label{Omega1unbounded}
    \{ z \in \mathbb C \mid  |z| > R, \, -\pi/2 + \varepsilon < \arg z < \pi/2 - \varepsilon \}
    & \subset \Omega_{+,\infty},
    \end{align}
    and
    \begin{align} \label{Omega2unbounded}
     \{ z \in \mathbb C \mid |z| > R, \, \pi/2 + \varepsilon < \arg z < 3 \pi/2 - \varepsilon \}
    & \subset \Omega_{-,\infty}. \end{align}
In particular,  $\Omega_{+,\infty}$ lies to the right of $\Omega_{-,\infty}$.
\item[\rm (b)]  Both $\Omega_{+,\infty}$ and $\Omega_{-,\infty}$ are connected.
\end{enumerate}
\end{lemma}

\begin{proof} Part (a) follows from \eqref{asymptotics:differencelambdas} and the
the fact that $a_k > a_{k+1}$. Part (b) follows from part (a) in a similar
way as in \cite[Proof of Lemma 2.4]{DelKui}, to which we refer for further
details. \end{proof}

\begin{lemma}\label{lemma:existencegloballenses}
We have
\begin{align} \label{alphaiinOmega}
    \alpha_i\in\Omega_{+,\infty} \qquad \text{and} \qquad
    \beta_{i+1}\in\Omega_{-,\infty},
\end{align} where we recall $i$ is related to $k$ as in
\eqref{verticalstep:condition1} and \eqref{verticalstep:condition2}.
\end{lemma}

\begin{proof}  By applying
the variational conditions \eqref{eulerlagrange:lambda2} twice, first with the
index $i$ and then with $i+1$, we obtain
\begin{equation} \label{ELtwice1}
    \Re(\lambda_{k+1}(x)-\lambda_k(x)) =
    \Re(\lambda_{k+1}(x)-\lambda_{p+l}(x)) \geq 0,\quad x\in [\alpha_i,\beta_i].
\end{equation}
and the inequality \eqref{ELtwice1} is strict for $x = \alpha_i$ because of
the statement about the strict inequality in Lemma~\ref{lem:eulerlagrangerestatement}.
Hence $\alpha_i\in\Omega_+$,
and in a similar way we obtain $\beta_{i+1}\in\Omega_-$.

To show that $\alpha_i$ belongs to the unbounded component of $\Omega_+$ we
argue as in \cite[Proof of Lemma 2.4]{DelKui}. What we use is that
$\bigcup_{j:\,  k(j)=k} [\alpha_j,\beta_j]$ lies to the right of $\bigcup_{j:\,
k(j)=k+1} [\alpha_j,\beta_j]$, and that $\alpha_i$ is the left-most point of
$\bigcup_{j:\,  k(j)=k} [\alpha_j,\beta_j]$, and that $\beta_{i+1}$ is the
right-most point of $\bigcup_{j:\,  k(j)=k+1} [\alpha_j,\beta_j]$.

The proof is by contradiction. Suppose that $\alpha_i$ does not belong
to the unbounded component of $\Omega_+$. Then the set
\begin{align}\label{def:omegaplusi}
    \Omega_{+,\alpha_i} &:=
        \{z\in\Omega_+ \mid \exists\textrm{ a connected path in $\Omega_+$ from $z$ to $\alpha_i$}\}
        \end{align}
is bounded, it is symmetric with respect to the real line, and it contains $\alpha_i$.
Also $\Re(\lambda_{k+1}-\lambda_k)$ is zero on the boundary of $\Omega_{+,\alpha_i}$ and
strictly positive inside of $\Omega_{+,\alpha_i}$ by construction. Since subharmonic
functions satisfy a maximum principle, it follows that $\Re(\lambda_{k+1}-\lambda_k)$
is not subharmonic on all of $\Omega_{+,\alpha_i}$. Then  by Lemma \ref{lemma:subsuperharmonisch}
we conclude that $\Omega_{+,\alpha_i}$ has a nonempty intersection with
$\bigcup_{j:\,  k(j)=k+1} [\alpha_j,\beta_j]$. Any point of intersection
lies strictly to the left of $\beta_{i+1}$, since $\beta_{i+1} \in \Omega_-$
and $\beta_{i+1}$ is the
right-most point of $\bigcup_{j:\,  k(j)=k+1} [\alpha_j,\beta_j]$.
Then, because of symmetry in the real axis, it follows that $\Omega_{+,\alpha_i}$
surrounds the point $\beta_{i+1}$. The set
\begin{align} \label{def:omegamini} \Omega_{-,\beta_{i+1}} &:= \{z\in\Omega_- \mid
    \exists\textrm{ a connected path in $\Omega_-$ from $z$ to $\beta_{i+1}$}\}
\end{align}
is then bounded, and it does not intersect with $\bigcup_{j:\,  k(j)=k} [\alpha_j,\beta_j]$.
By Lemma \ref{lemma:subsuperharmonisch}, we then have that $\Re(\lambda_{k+1}-\lambda_k)$ is superharmonic on
$\Omega_{-,\beta_{i+1}}$.
However, $\Re(\lambda_{k+1}-\lambda_k)$ is zero on the boundary, and strictly negative inside
of $\Omega_{-,\beta_{i+1}}$, which gives a contradiction with the minimum principle for
superharmonic functions.

Thus  $\alpha_i$ belongs
to the unbounded component of $\Omega_+$, and similarly $\beta_{i+1}$ belongs to
the unbounded component of $\Omega_-$.
 \end{proof}

Let $X_0>0$ be such that
\[ \Re(\lambda_{k'}(x)-\lambda_{p+l'}(x)) > 0, \]
for all $x\in(-\infty,-X_0) \cup (X_0,\infty)$ and all $k' = 1, \ldots, p$,
$l' = 1, \ldots, q$.
The existence of such a constant $X_0$
follows from \eqref{asymptotics:lambdaa}--\eqref{asymptotics:lambdab}.

The next result is an immediate consequence of Lemma
\ref{lemma:existencegloballenses}.

\begin{theorem}\label{theorem:globallensesexist}
There exist two simple closed contours
\[ \Gamma_{+,i}\subset\Omega_{+,\infty} \qquad  \textrm{and} \qquad
    \Gamma_{-,i}\subset\Omega_{-,\infty} \]
    such that
\begin{itemize}
\item[\rm (a)] $\Gamma_{+,i}$ surrounds the interval $[\alpha_i,\beta_1]$.
\item[\rm (b)] $\Gamma_{-,i}$ surrounds the interval $[\alpha_{p+q-1}, \beta_{i+1}]$.
\item[\rm (c)] both $\Gamma_{-,i}$ and $\Gamma_{+,i}$ intersect the interval
$(\beta_{i+1},\alpha_{i})$ in exactly one point, which we denote by $x_i$,
$y_i$, respectively. We have $x_i<y_i$.
\item[\rm (d)] both $\Gamma_{-,i}$ and $\Gamma_{+,i}$ have one extra intersection point with
the real axis, which lies inside the interval $(-\infty,-X_0)$, $(X_0,\infty)$,
respectively.
\end{itemize}
As in \cite{ABK}, the contours $\Gamma_{+,i}$ and $\Gamma_{-,i}$ will be called
\emph{global lenses}.
\end{theorem}

\begin{remark} Instead of taking closed contours, one might also take
$\Gamma_{+,i}$ and $\Gamma_{-,i}$ to be unbounded, tending to infinity in the
right and left half of the complex plane, respectively, and both intersecting
the real line in exactly one point in the line segment
$(\beta_{i+1},\alpha_i)$. This is the construction that was used in
\cite{DelKui}.
\end{remark}

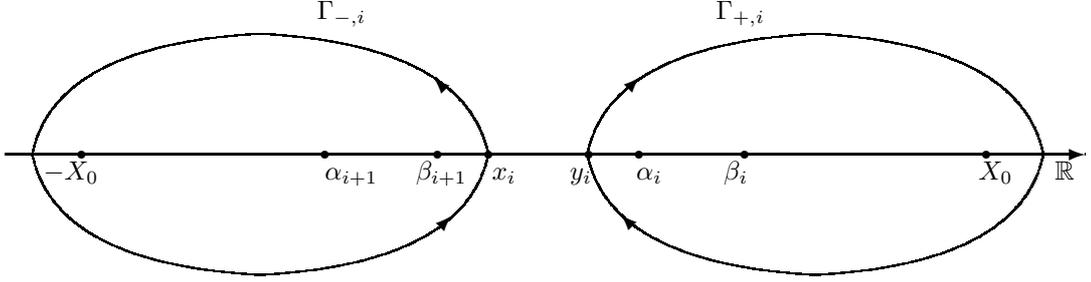
\begin{figure}[t]
\begin{center}
   \setlength{\unitlength}{1truemm}
   \begin{picture}(100,70)(-5,2)

       \qbezier(54.5,40)(57.5,25.5)(84.5,24)
       \qbezier(84.5,24)(111.5,25.5)(114.5,40)
       \qbezier(54.5,40)(57.5,54.5)(84.5,56)
       \qbezier(84.5,56)(111.5,54.5)(114.5,40)
       \put(71.5,58){$\Gamma_{+,i}$}

       \put(60,31){\thicklines\vector(-1,1){1}}
       \put(60.4,49.2){\thicklines\vector(1,1){1}}

       \qbezier(-18.5,40)(-15.5,25.5)(11.5,24)
       \qbezier(11.5,24)(38.5,25.5)(41.5,40)
       \qbezier(-18.5,40)(-15.5,54.5)(11.5,56)
       \qbezier(11.5,56)(38.5,54.5)(41.5,40)
       \put(19,58){$\Gamma_{-,i}$}

       \put(35.7,30.8){\thicklines\vector(1,1){1}}
       \put(35.4,49.1){\thicklines\vector(-1,1){1}}

       \put(-22,40){\line(1,0){140}}
       \put(116,36.6){$\mathbb R $}
       \put(120,40){\thicklines\vector(1,0){.0001}}
       \put(-12,40){\thicklines\circle*{1}}
       \put(-17,36.6){$-X_{0}$}
       \put(20.0,40){\thicklines\circle*{1}}
       \put(20,36.6){$\alpha_{i+1}$}
       \put(34.8,40){\thicklines\circle*{1}}
       \put(32,36.6){$\beta_{i+1}$}
       \put(41.5,40){\thicklines\circle*{1}}
       \put(42,36.6){$x_i$}
       \put(54.7,40){\thicklines\circle*{1}}
       \put(52.2,36.6){$y_i$}
       \put(61.3,40){\thicklines\circle*{1}}
       \put(61,36.6){$\alpha_i$}
       \put(75.2,40){\thicklines\circle*{1}}
       \put(72.5,36.6){$\beta_i$}
       \put(107,40){\thicklines\circle*{1}}
       \put(106,36.6){$X_{0}$}

   \end{picture}
   \vspace{-20mm}
   \caption{The figure shows how to open global lenses between the two intervals
   $[\alpha_{i+1},\beta_{i+1}]$ (left) and $[\alpha_i,\beta_i]$ (right).}
   \label{fig:globallenses}
\end{center}
\end{figure}

An illustration of Theorem \ref{theorem:globallensesexist} is shown in
Figure~\ref{fig:globallenses}. When convenient we also define $x_{0}=X_0$,
$y_0=\infty$, $x_{p+q-1}=-\infty$ and $y_{p+q-1}=-X_0$.

\subsection{Construction of global lenses: Horizontal step}
\label{subsection:horizontalstep}

Next we construct the global lenses for a horizontal step along the lattice
path. Thus assume that $i\in\{1,\ldots,p+q-2\}$ is such that
\begin{align}
\label{horizontalstep:condition1}
(k(i),l(i)) &= (k,l),\\
\label{horizontalstep:condition2} (k(i+1),l(i+1)) &= (k,l+1),
\end{align}
for certain $k=1,\ldots,p$ and $l=1,\ldots,q-1$.

The analysis of Section \ref{subsection:verticalstep} can be adapted to the
present case. Let us discuss the main points. We now define the open sets
$\Omega_+$ and $\Omega_-$ as follows
\begin{align}\label{def:omegaplusbis}
    \Omega_+ & := \{z\in\mathbb C\mid  \Re(\lambda_{p+l}(z)-\lambda_{p+l+1}(z))>0\} \\
    \label{def:omegaminbis} \Omega_- & := \{z\in\mathbb C\mid
    \Re(\lambda_{p+l}(z)-\lambda_{p+l+1}(z))<0\}.
\end{align}
The unbounded regions $\Omega_{+,\infty},\Omega_{-,\infty}\subset\mathbb C$ are
again defined as in \eqref{def:omegaplusinfty}--\eqref{def:omegamininfty}. Then
Lemma \ref{lemma:infinitecomponents} remains valid.

Using the above definitions, Lemma~\ref{lemma:existencegloballenses} and
Theorem~\ref{theorem:globallensesexist} both remain valid as well. Thus we can
construct the global lenses $\Gamma_{+,i}$ and $\Gamma_{-,i}$ in exactly the
same way as before.

\subsection{Gaussian elimination step}
\label{subsection:globallensesalgorithm}

Now we will show how to apply Gaussian elimination inside the global lenses.
First, it might be worthwhile to recall the basic idea of Gaussian elimination
in our present context, see e.g.\ \cite{GVL}. Let
\begin{equation}\label{rankonematrix:general} J=\mathbf{u}\mathbf{v}^T =
\begin{pmatrix}u_1 \\ \vdots \\ u_p\end{pmatrix}
\begin{pmatrix}v_1 & \ldots & v_q\end{pmatrix}
\end{equation}
be a rank-one matrix of size $p$ by $q$. Assume that one multiplies $J$ on the
left with the matrix \begin{equation}\label{gausselim:basicidea1}
I_p-\frac{u_{k+1}}{u_{k}}\vece_{k+1}\vece^T_{k}.\end{equation} Here and below
we use $\vec{e}_k$ to denote the column vector with all entries equal to
zero, except for the $k$th entry which equals 1. The length of $\vece_k$ will
be clear from the context. Note that the outer product
$\vece_{k+1}\vece_k^T$ in \eqref{gausselim:basicidea1} is the matrix with all zero entries
except for the $(k+1,k)$ entry which equals one.

The multiplication with \eqref{gausselim:basicidea1} on the left is equivalent
to applying an elementary row operation to the rows of $J$, where to row $k+1$
one adds $-\frac{u_{k+1}}{u_{k}}$ times row $k$. This row operation is such
that the entries of row $k+1$ of $J$ are eliminated.

Similarly, assume that one multiplies $J$ on the right with the matrix
\begin{equation}\label{gausselim:basicidea2} I_q-\frac{v_{l+1}}{v_{l}}\vece_{l}\vece^T_{l+1}.
\end{equation} This is equivalent to
applying an elementary column operation to the columns of $J$, where to column
$l+1$ one adds $-\frac{v_{l+1}}{v_{l}}$ times column $l$. This column operation
is such that the entries of column $l+1$ of $J$ are eliminated.

Let us see how we can apply these ideas in the present context. The role of the
rank-one matrix \eqref{rankonematrix:general} of size $p$ by $q$ will be played
by the top right submatrix $J_{1,2}(x)$ of the jump matrix \eqref{jumpsX1},
cf.\ \eqref{jumpsX4}. The mechanism to multiply the jump matrix on the left or
on the right with a transformation matrix of the form
\eqref{gausselim:basicidea1}--\eqref{gausselim:basicidea2} is to define a new
RH matrix $T(z) = X(z)B(z)$ for a suitable transformation matrix $B(z)$ and
subsequently apply Lemma \ref{lemma:transformation:jumps}.

Now we are ready to describe the Gaussian elimination in detail. This will be
the next transformation $X\mapsto T$ in the steepest descent analysis of the RH
problem.

\begin{algorithm}\label{algorithm:twosweeps} (The transformation $X \mapsto T$)
\begin{enumerate}
\item (Initialization.)
We initialize $T(z):=X(z)$. \item (Forward sweep.) For each $i=1,\ldots,p+q-2$
we open the global lens $\Gamma_{-,i}$ in Theorem
\ref{theorem:globallensesexist} and we update, in case of a vertical step
\eqref{verticalstep:condition1}--\eqref{verticalstep:condition2},
$$
T(z) = \left\{\begin{array}{ll}
T(z)\left(I_{p+q}+\exp(n(\lambda_{k+1}(z)-\lambda_{k}(z)))\vece_{k}\vece^T_{k+1}\right),&
\textrm{ inside the lens }\Gamma_{-,i}\\
T(z),& \textrm{ elsewhere},\end{array}\right.
$$
and in case of a horizontal step
\eqref{horizontalstep:condition1}--\eqref{horizontalstep:condition2},
$$
T(z) = \left\{\begin{array}{ll}
T(z)\left(I_{p+q}-\exp(n(\lambda_{p+l}(z)-\lambda_{p+l+1}(z)))\vece_{p+l+1}\vece^T_{p+l}\right),& \textrm{ inside }\Gamma_{-,i}\\
T(z),& \textrm{ elsewhere}.\end{array}\right.
$$
\item (Backward sweep.)
For each $i=p+q-2,\ldots,2,1$ we open the global lens $\Gamma_{+,i}$ in Theorem
\ref{theorem:globallensesexist} and we update, in case of a vertical step
\eqref{verticalstep:condition1}--\eqref{verticalstep:condition2},
$$
T(z) = \left\{\begin{array}{ll}
T(z)\left(I_{p+q}+\exp(n(\lambda_{k}(z)-\lambda_{k+1}(z)))\vece_{k+1}\vece^T_{k}\right),&
\textrm{ inside the lens }\Gamma_{+,i}\\
T(z),& \textrm{ elsewhere},\end{array}\right.
$$
and in case of a horizontal step
\eqref{horizontalstep:condition1}--\eqref{horizontalstep:condition2},
$$
T(z) = \left\{\begin{array}{ll}
T(z)\left(I_{p+q}-\exp(n(\lambda_{p+l+1}(z)-\lambda_{p+l}(z)))\vece_{p+l}\vece^T_{p+l+1}\right),& \textrm{ inside }\Gamma_{+,i}\\
T(z),& \textrm{ elsewhere}.\end{array}\right.
$$
\end{enumerate}
\end{algorithm}

Incidently, we note that the forward and backward sweeps in the above algorithm
commute. But one is not allowed to change the order in which the index $i$
varies inside the sweeps.

It is easy to see that Algorithm \ref{algorithm:twosweeps} does not change the
jump matrix in \eqref{jumpsX1} except for its top right submatrix $J_{1,2}(x)$.
To see what happens with the latter, we have to distinguish between different
regions of the complex plane.

First assume that $x$ belongs to one of the intervals $(y_{j},x_{j-1})\supset
[\alpha_j,\beta_j]$, $j=1,\ldots,p+q-1$. (Recall that $y_{p+q-1}=-X_0$ and $x_{0}=X_0$.)
From Theorem \ref{theorem:globallensesexist} we see that $x$ lies inside the global
lens $\Gamma_{-,i}$ precisely when $i=1,\ldots,j-1$. Hence the \lq relevant\rq\
indices in the forward sweep in Algorithm~\ref{algorithm:twosweeps} are
$i=1,\ldots,j-1$. During the corresponding operations, the entries in rows
$1,2,\ldots,k(j)-1$ and columns $1,2,\ldots,l(j)-1$ of the matrix $J_{1,2}(x)$
are cancelled by Gaussian elimination.

On the other hand, Theorem \ref{theorem:globallensesexist} shows that $x$ lies
inside the global lens $\Gamma_{+,i}$ precisely when $i=p+q-2,p+q-3,\ldots,j$.
Hence the relevant indices in the backward sweep in
Algorithm~\ref{algorithm:twosweeps} are $i=p+q-2,p+q-3,\ldots,j$. During the
corresponding operations, the entries in rows $p,p-1,\ldots,k(j)+1$ and columns
$q,q-1,\ldots,l(j)+1$ of $J_{1,2}(x)$ are cancelled by Gaussian elimination.

It follows that at the end of the two sweeps in
Algorithm~\ref{algorithm:twosweeps}, all the entries of the rank-one matrix
$J_{1,2}(x)$ are eliminated, except for the $(k(j),l(j))$ entry which equals
$$\exp(n(\lambda_{p+l(j),+}(x)-\lambda_{k(j),-}(x))).$$

Recall that in the above description, we assumed that $x\in
(y_{j},x_{j-1})\supset [a_j,b_j]$. Next, let us assume that $x$ belongs to one of the
gaps $(x_{j},y_{j})$ for certain $j$. We can then repeat the above arguments
and find that at the end of Algorithm~\ref{algorithm:twosweeps}, all the
entries of $J_{1,2}(x)$ are eliminated except for two of them. In case of a
vertical step \eqref{verticalstep:condition1}--\eqref{verticalstep:condition2}
these are the $(k(j),l(j))$ and $(k(j)+1,l(j))$ entries, which are given by
$$
\exp(n(\lambda_{p+l(j)}(x)-\lambda_{k(j)}(x))),\quad
\exp(n(\lambda_{p+l(j)}(x)-\lambda_{k(j)+1}(x))),
$$
respectively. But by the variational inequality in
\eqref{eulerlagrange:lambda2}, which is strict according to
Lemma~\ref{lemma:eulerlagrange:strictinequality}, we see that both entries are
exponentially small for $n\to\infty$. A similar argument applies in case of a
horizontal step
\eqref{horizontalstep:condition1}--\eqref{horizontalstep:condition2}.

Finally we should note that by the operations in
Algorithm~\ref{algorithm:twosweeps}, the RH matrix $T(z)$ also has a jump on
each of the contours $\Gamma_{+,i}$ and $\Gamma_{-,i}$. For example, in case
of a vertical step
\eqref{verticalstep:condition1}--\eqref{verticalstep:condition2} the jump
matrix on the contour $\Gamma_{+,i}$ takes the form
\[ I_{p+q}\pm \exp(n(\lambda_{k}(z)-\lambda_{k+1}(z)))\vece_{k+1}\vece^T_{k},\]
see Algorithm \ref{algorithm:twosweeps}.
But by our assumption that $\Gamma_{+,i}\subset\Omega_+$ we see that this jump
matrix is uniformly exponentially close to the identity matrix when
$n\to\infty$. A similar argument holds for the jumps along the contours
$\Gamma_{-,i}\subset\Omega_-$.

Summarizing, we established that  $T$ satisfies the
following RH problem.
\begin{rhp}\label{RHP:T}\textrm{ }
\begin{itemize}
\item[\rm (1)] $T$ is analytic
on $\mathbb C \setminus(\mathbb R\cup\bigcup_{i=1}^{p+q-2}(\Gamma_{+,i}\cup
\Gamma_{-,i}))$.
 \item[\rm (2)] For $x \in \mathbb R \cup \bigcup_{i=1}^{p+q-2}(\Gamma_{+,i}\cup \Gamma_{-,i}))$ we have that
\begin{equation}\label{jumpsT1}
T_{+}(x) = T_{-}(x) J_T(x)
\end{equation}
where $J_T(x)$ satisfies the following
\begin{enumerate}
\item[\rm (a)]
 For $x\in (y_i,x_{i-1})\supset [\alpha_i,\beta_i]$, $i= 1,\ldots,p+q-1$,
 we have that $J_T(x)$ equals the identity matrix, except for the $2\times 2$ block
lying on the intersection of rows and columns $k(i)$ and $p + l(i)$, which is
given by
\begin{equation}\label{jumpsT3}
\begin{pmatrix}
\exp(n(\lambda_{k(i),+}(x)-\lambda_{k(i),-}(x))) &
\exp(n(\lambda_{p+l(i),+}(x)-\lambda_{k(i),-}(x))) \\
0& \exp(n(\lambda_{p+l(i),+}(x)-\lambda_{p+l(i),-}(x)))
\end{pmatrix}.
\end{equation}
\item[\rm (b)] For $x \in (-\infty, y_{p+q-1}) \cup \left( \bigcup_i (x_i,y_i) \right) \cup (x_0, \infty)$
we have that $J_T(x)$ is  exponentially close
to the identity matrix as $n\to\infty$, both uniformly as well as in $L^2$ sense.
\item[\rm (c)]
For $x\in\bigcup_{i}(\Gamma_{+,i}\cup \Gamma_{-,i})$, the jump matrix $J_T$ is also
exponentially close to the identity matrix as $n\to\infty$ in uniform sense (and
therefore also in $L^2$ since the contours $\Gamma_{\pm,i}$ are compact).
\end{enumerate}
\item[\rm (3)] As $z\to\infty$, we have that
\begin{equation}\label{asymptoticconditionT}
    T(z) =I_{p+q}+O(1/z).
\end{equation}
\end{itemize}
\end{rhp}

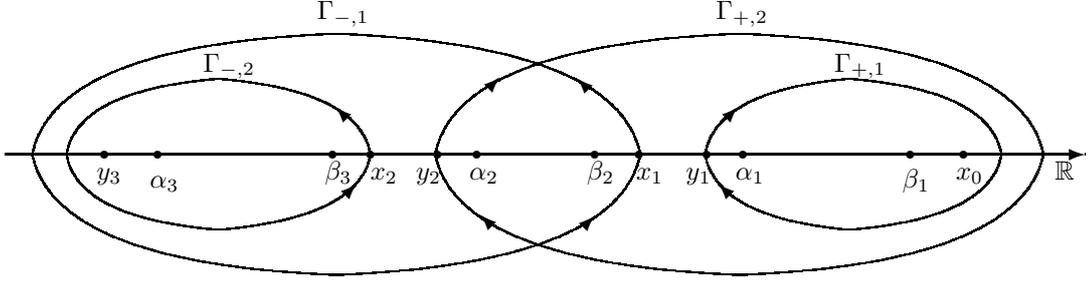
\begin{figure}[t]
\begin{center}
   \setlength{\unitlength}{1truemm}
   \begin{picture}(100,70)(-5,2)
       \qbezier(-14,40)(-12.5,31.5)(6,30)
       \qbezier(6,30)(24.5,31.5)(26,40)
       \qbezier(-14,40)(-12.5,48.5)(6,50)
       \qbezier(6,50)(24.5,48.5)(26,40)
       \put(4,51){$\Gamma_{-,2}$}
       \qbezier(34.5,40)(37.5,25.5)(74.5,24)
       \qbezier(74.5,24)(111.5,25.5)(114.5,40)
       \qbezier(34.5,40)(37.5,54.5)(74.5,56)
       \qbezier(74.5,56)(111.5,54.5)(114.5,40)
       \put(71.5,58){$\Gamma_{+,2}$}

       \put(23.2,35){\thicklines\vector(1,1){1}}
       \put(22.5,45.2){\thicklines\vector(-1,1){1}}
       \put(41.5,31){\thicklines\vector(-1,1){1}}
       \put(42,49.2){\thicklines\vector(1,1){1}}

       \qbezier(-18.5,40)(-15.5,25.5)(21.5,24)
       \qbezier(21.5,24)(58.5,25.5)(61.5,40)
       \qbezier(-18.5,40)(-15.5,54.5)(21.5,56)
       \qbezier(21.5,56)(58.5,54.5)(61.5,40)
       \put(19,58){$\Gamma_{-,1}$}
       \qbezier(70,40)(71.5,31.5)(89,30)
       \qbezier(89,30)(107.5,31.5)(109,40)
       \qbezier(70,40)(71.5,48.5)(89,50)
       \qbezier(89,50)(107.5,48.5)(109,40)
       \put(87,51){$\Gamma_{+,1}$}

       \put(72.8,35){\thicklines\vector(-1,1){1}}
       \put(73.3,45.2){\thicklines\vector(1,1){1}}
       \put(54.5,30.8){\thicklines\vector(1,1){1}}
       \put(54.2,49.1){\thicklines\vector(-1,1){1}}

       \put(-22,40){\line(1,0){140}}
       \put(116,36.6){$\mathbb R $}
       \put(120,40){\thicklines\vector(1,0){.0001}}
       \put(-9,40){\thicklines\circle*{1}}
       \put(-10,36.6){$y_{3}$}

       \put(-2,40){\thicklines\circle*{1}}
       \put(-3,35.6){$\alpha_{3}$}
       \put(21,40){\thicklines\circle*{1}}
       \put(20,36.6){$\beta_{3}$}
       \put(26.0,40){\thicklines\circle*{1}}
       \put(26,36.6){$x_2$}
       \put(34.8,40){\thicklines\circle*{1}}
       \put(32,36.6){$y_2$}
       \put(40,40){\thicklines\circle*{1}}
       \put(39,36.6){$\alpha_2$}
       \put(55.5,40){\thicklines\circle*{1}}
       \put(54.5,36.6){$\beta_2$}
       \put(61.3,40){\thicklines\circle*{1}}
       \put(61,36.6){$x_1$}
       \put(70.2,40){\thicklines\circle*{1}}
       \put(67.5,36.6){$y_1$}
       \put(75,40){\thicklines\circle*{1}}
       \put(74,36.6){$\alpha_1$}
       \put(97,40){\thicklines\circle*{1}}
       \put(96,35.6){$\beta_1$}
       \put(104,40){\thicklines\circle*{1}}
       \put(103,36.6){$x_{0}$}

   \end{picture}
 \vspace{-20mm}
   \caption{The figure shows the contours in the RH problem for the matrix $T$ for the example
   of Figures  \ref{fig:graph1} and \ref{fig:Riemannsurface4x4}.
   We have three intervals $[\alpha_i,\beta_i]$, $i=1,2,3$, and four global
   lenses $\Gamma_{\pm,i}$, $i=1,2$.}
   \label{fig:globallenses:example}
\end{center}
\end{figure}

Let us illustrate this RH problem for the example in Figures \ref{fig:graph1}
and \ref{fig:Riemannsurface4x4}. Then we have three intervals
$[\alpha_i,\beta_i]$, $i=1,2,3$, and four global lenses between them: see
Figure \ref{fig:globallenses:example}.

The jump conditions \eqref{jumpsT1}--\eqref{jumpsT3} can now be written as
\begin{equation*}
T_{+} = T_{-}
\begin{pmatrix}
\exp(n(\lambda_{1,+}-\lambda_{1,-})) & 0 & \exp(n(\lambda_{3,+}-\lambda_{1,-})) &0 \\
0 & 1 & 0 & 0 \\
0 & 0 & \exp(n(\lambda_{3,+}-\lambda_{3,-})) &0 \\
0 & 0 & 0 & 1
\end{pmatrix}
\end{equation*}
on the interval $(y_1,x_0)\supset [\alpha_1, \beta_1]$,
\begin{equation*}
T_{+} = T_{-}
\begin{pmatrix}
1 & 0 & 0 & 0 \\
0 & \exp(n(\lambda_{2,+}-\lambda_{2,-})) & \exp(n(\lambda_{3,+}-\lambda_{2,-})) & 0\\
0 & 0 & \exp(n(\lambda_{3,+}-\lambda_{3,-})) & 0 \\
0 & 0 & 0 & 1
\end{pmatrix}
\end{equation*}
on the interval $(y_2,x_1) \supset [\alpha_2,\beta_2]$, and
\begin{equation*}
T_{+} = T_{-}
\begin{pmatrix}
1 & 0 & 0 & 0 \\
0 & \exp(n(\lambda_{2,+}-\lambda_{2,-})) & 0 & \exp(n(\lambda_{4,+}-\lambda_{2,-})) \\
0 & 0 & 1 & 0 \\
0 & 0 & 0 & \exp(n(\lambda_{4,+}-\lambda_{4,-})) \\
\end{pmatrix}
\end{equation*}
on the interval $(y_3,x_2)\supset [\alpha_3,\beta_3]$. The jump matrices on the remaining
contours $(-\infty,y_3)$, $(x_2,y_2)$,  $(x_1,y_1)$, $(x_0,\infty)$,
$\Gamma_{+,1}$, $\Gamma_{-,1}$, $\Gamma_{+,2}$ and $\Gamma_{-,2}$ in Figure
\ref{fig:globallenses:example} are all exponentially close to the
identity matrix as $n\to\infty$.

\section{Final transformations of the RH problem}
\label{section:remainingstepssteepestdescent}

The jump matrices $J_T$ in the RH problem for $T$ are nontrivial
only in the $2 \times 2$ block given by \eqref{jumpsT3} on
the interval $[\alpha_i,\beta_i]$.

\subsection{Third transformation: Opening of the local lenses}
\label{subsection:locallenses}

In the transformation $T\mapsto S$ of the RH problem we transform the
oscillatory entries of the jump matrix \eqref{jumpsT1}--\eqref{jumpsT3} along
each interval $[\alpha_i,\beta_i]$ into exponentially decaying ones. To this end we open a
local lens around the interval $[\alpha_i,\beta_i] = [\alpha_i,\beta_i]$. Since the RH problem
is locally of size $2$ by $2$ this can be done in the standard way \cite{Dei,DKMVZ1}.

\begin{figure}[t]
\begin{center}
   \setlength{\unitlength}{1truemm}
   \begin{picture}(100,70)(-5,2)

       \put(0,40){\line(1,0){95}}
       \put(93,36.6){$\mathbb R $}
       \put(95,40){\thicklines\vector(1,0){.0001}}
       \put(28,40){\thicklines\circle*{1}}
       \put(5,36.6){$y_{i}$}
       \put(6,40){\thicklines\circle*{1}}
       \put(27,36.6){$\alpha_i$}
       \put(62,40){\thicklines\circle*{1}}
       \put(61,36.6){$\beta_i$}
       \put(86,40){\thicklines\circle*{1}}
       \put(85,36.6){$x_{i-1}$}

       \qbezier(28,40)(45,60)(62,40)
       \put(45,50){\thicklines\vector(1,0){1}}
       \put(45,52){$L_{+,i}$}

       \qbezier(28,40)(45,20)(62,40)
       \put(45,30){\thicklines\vector(1,0){1}}
       \put(45,26){$L_{-,i}$}
   \end{picture}
 \vspace{-20mm}
   \caption{A local lens around the interval $[\alpha_i,\beta_i]$.}
   \label{fig:locallenses}
\end{center}
\end{figure}
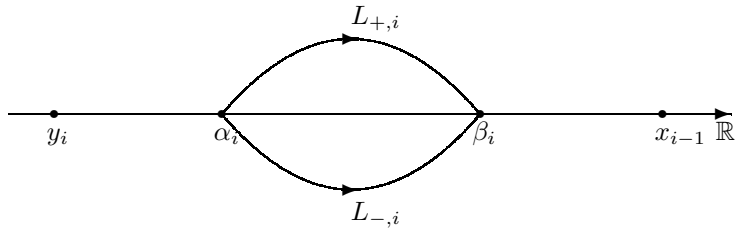

For each $i = 1,\ldots, p+q-1$ we open a lens around the interval
$[\alpha_i,\beta_i]$, $i=1,\ldots,p+q-1$ as in Figure \ref{fig:locallenses}
so that
\[ \Re (\lambda_{k(i)} - \lambda_{p+l(i)} ) < 0 \]
on the lips of the lens. It follows from an argument based
on the Cauchy-Riemann equations that this is indeed possible, cf. \cite{Dei}.
The lenses are small (so called local lenses),
and in particular they do not intersect with the global lenses $\Gamma_{\pm,i}$.
We use $L_{+,i}$ and $L_{-,i}$ to denote the upper and lower lip of the lens,
respectively.

We define the matrix function $S$  as follows
\begin{equation}\label{def:S}
S(z) = \left\{\begin{array}{l}
    T(z)\left(I_{p+q}-\exp(n(\lambda_{k(i)}(z)-\lambda_{p+l(i)}(z)))\vece_{p+l(i)}\vece^T_{k(i)}\right), \\
    \qquad \qquad \textrm{ in the upper part of the lens around $[\alpha_i,\beta_i]$,} \\
T(z)\left(I_{p+q}+\exp(n(\lambda_{k(i)}(z)-\lambda_{p+l(i)}(z)))\vece_{p+l(i)}\vece^T_{k(i)}\right), \\
    \qquad  \qquad\textrm{ in the lower part of the lens around $[\alpha_i,\beta_i]$,}\\
T(z), \qquad \textrm{outside of all the lenses}.\end{array}\right.
\end{equation}

Then $S$ satisfies the following RH problem.

\begin{rhp}\label{RHP:S}\textrm{ }
\begin{itemize}
\item[\rm (1)] $S$ is analytic
in $\mathbb C \setminus\left(\mathbb R\cup\bigcup_{i}
(\Gamma_{+,i}\cup\Gamma_{-,i})\cup \bigcup_{i} (L_{+,i}\cup L_{-,i} )\right)$;
 \item[\rm (2)] For $x \in \mathbb R \cup\bigcup_{i}
(\Gamma_{+,i}\cup\Gamma_{-,i})\cup \bigcup_{i} (L_{+,i}\cup L_{-,i} )$ we have
that
\begin{equation}\label{jumpsS1}
S_{+}(x) = S_{-}(x) J_S(x)
\end{equation}
where $J_S(x)$ satisfies the following
\begin{enumerate}
\item[\rm (a)]
For $x \in [y_i, x_{i-1}] \cup L_{+,i} \cup L_{-,i}$ we have that
$J_S(x)$ is the identity matrix except for the $2\times 2$ block
on the intersection of rows and columns $k(i)$ and $p + l(i)$, which is
given by
\begin{equation}\label{jumpsS3}
\begin{pmatrix} 0 & 1 \\ -1 & 0\end{pmatrix}, \qquad \text{for } x\in [\alpha_i,\beta_i],
\end{equation}
\begin{equation} \label{jumpsS4}
\begin{pmatrix} 1 & 0 \\ \exp(n(\lambda_{k(i)}-\lambda_{p+l(i)})) & 1\end{pmatrix},
    \qquad \text{for }  x\in L_{\pm,i},
    \end{equation}
and by
\begin{equation} \label{jumpsS5}
\begin{pmatrix} \exp(n(\lambda_{k(i),+}(x)-\lambda_{k(i),-}(x)))
    & \exp(n(\lambda_{p+l(i),+}(x)-\lambda_{k(i),-}(x))) \\
    0& \exp(n(\lambda_{p+l(i),+}(x)-\lambda_{p+l(i),-}(x)))
\end{pmatrix},
\end{equation}
for $ x \in [y_i,\alpha_i] \cup [\beta_i, x_{i-1}]$.
\item[\rm (b)]
On the other parts of the contour we have $J_S(x) = J_T(x)$ and $J_S(x)$ is
exponentially close to the identity matrix as $n \to \infty$, both uniformly
and in the $L^2$ sense.
\end{enumerate}
\item[\rm (3)] As $z\to\infty$, we have that
\begin{equation}\label{asymptoticconditionS}
    S(z) =I_{p+q}+O(1/z).
\end{equation}
\end{itemize}
\end{rhp}

From standard arguments based on the Cauchy-Riemann conditions \cite{Dei} it
follows that the local lenses $L_i$ can be chosen so that the jumps on $L_i$ in
\eqref{jumpsS4} are uniformly exponentially close to the identity matrix, away
for a neighborhood of the endpoints $\alpha_i$, $\beta_i$.

\subsection{Global parametrix}
\label{subsection:globalparametrix}

In this subsection we build a global parametrix $P^{(\infty)}(z)$, which will be
a good approximation to the RH problem away from the endpoints
$\alpha_i$, $\beta_i$, $i=1,\ldots,p+q-1$. The construction will be quite similar
to the one in \cite{DKV}.

We will construct the matrix function $P^{(\infty)}(z)$ such that it satisfies
the following RH problem, obtained from RH problem \ref{RHP:S} by ignoring all
exponentially small entries of the jump matrices.

\begin{rhp}\label{RHP:Pinfty}\textrm{ }
\begin{itemize}
\item[(1)] $P^{(\infty)}(z)$ is analytic
in $\mathbb C \setminus\bigcup_{i=1}^{p+q-1} [\alpha_i,\beta_i] $;
 \item[(2)] For $x\in \bigcup_{i} (\alpha_i,\beta_i)$,  we have that
\begin{equation}\label{jumpsPinfty1}
P^{(\infty)}_{+}(x) = P^{(\infty)}_{-}(x) J_{P^{(\infty)}}(x)
\end{equation}
where the jump matrix $J_{P^{(\infty)}}(x)$ equals
\begin{equation}\label{jumpsPinfty2}
J_{P^{(\infty)}}(x) =
\begin{pmatrix}
I_{k(i)-1} & 0 & 0 & 0 &0 \\
0& 0 & 0& 1 &0 \\
0&0 & I_{p+l(i)-1-k(i)} & 0&0 \\
0& -1 &0 & 0 &0 \\
0&0 &0 &0 & I_{q-l(i)}
\end{pmatrix},
\end{equation}
for $x \in (\alpha_i,\beta_i)$.
\item[(3)] As $z\to\infty$, we have that
\begin{equation*}\label{asymptoticconditionPinfty}
    P^{(\infty)}(z) = I_{p+q}+O(1/z).
\end{equation*}
\end{itemize}
\end{rhp}

To solve this RH problem, we will use the fact that the Riemann surface $\mathcal R$
in Section \ref{subsection:xifunctions} has genus zero. From general algebraic
geometry \cite{Mir}, this implies the existence of a rational parametrization
\begin{align}\label{ratpar:xi:z}
    \xi=\xi(v), \qquad z=z(v) \end{align}
where $v$ runs through the extended complex plane $\overline{\mathbb C}$ (Riemann
sphere).

The $v$-plane is then partitioned into $p+q$ disjoint open sets $\Omega_j$,
$j=1,\ldots,p+q$, where $\Omega_j$ is defined as the inverse image under
\eqref{ratpar:xi:z} of the $j$th sheet $\mathcal R_j$ of the Riemann surface.
Correspondingly we have $p+q$ inverse functions $v_j(z)$ of \eqref{ratpar:xi:z}
such that
\begin{equation}\label{defvmappings} v_j: \mathcal R_j\to\Omega_j,\quad
j=1,\ldots,p+q,
\end{equation}
is a bijection.
We use  $v_j(\infty)$ to denote the image under this map of the point at infinity of
the $j$th sheet $\mathcal R_j$, $j=1,\ldots,p+q$. Hence $v_j(\infty)\in\Omega_j$.

For $i = 1, \ldots, p+q-1$, the common boundary of $\Omega_{k(i)}$ and $\Omega_{p+l(i)}$ in the
$v$-plane, is an analytic curve $\ccal_{i}$ with a natural partition
\begin{equation}
    \ccal_{i} = \ccal_{+,i}\cup \ccal_{-,i},\end{equation}
where $\ccal_{+,i}$  is the image of the upper side of the cut
$[\alpha_i,\beta_i]$ under the mapping $v_{k_i}$, and $\ccal_{-,i}$ is the
image of the lower side. The two parts $\ccal_{\pm,i}$ meet at two points
$\gamma^{(1)}_{i}$ and $\gamma^{(2)}_{i}$, which are the images of the
endpoints $\alpha_i$ and $\beta_i$ of the cut $[\alpha_i,\beta_i]$,
respectively.

Define the polynomial
\begin{equation}\label{def:g}
    g(v)=\prod_{i=1}^{p+q-1}\left((v-\gamma^{(1)}_{i})(v-\gamma^{(2)}_{i})\right),
\end{equation}
and its square root
\[ \sqrt{g(v)} \]
which is defined as an analytic function in
the $v$-plane, with a cut along the disjoint union of arcs
$\bigcup_{i}\ccal_{+,i}$. We assume $\sqrt{g(v)} \sim v^{p+q-1}$ as $v \to \infty$.

We then construct a global parametrix $P^{(\infty)}(z)$ as in \cite{DKV}.
We define for $z \in \mathbb C \setminus \bigcup_i [\alpha_i, \beta_i]$,
\begin{equation}\label{globalparametrix}
P^{(\infty)}(z) =
\begin{pmatrix} \ds f_i(v_j(z)) \end{pmatrix}_{i,j=1}^{p+q},
    \qquad f_i(v) = \frac{l_i(v)}{\sqrt{g(v)}},
\end{equation}
where  $l_i$ is the Lagrange interpolation polynomial for the
points $v_1(\infty), \ldots v_{p+q}(\infty)$. That is, $l_i$ is a polynomial
of degree $p+q-1$ so that
\[ f_i(v_j(\infty)) = \delta_{i,j}, \qquad j=1,\ldots,p+q. \]

The fact that $P^{(\infty)}(z)$ in \eqref{globalparametrix} satisfies
conditions (1) and (3) in RH problem~\ref{RHP:Pinfty} is immediate. For the
jump condition (2) we need to show that
\begin{align}
\left. \begin{array}{l} f_i(v_{k(i),+}(x)) = -f_i(v_{p+l(i),-}(x)) \\
    f_i(v_{p+l(i),+}(x)) = f_i(v_{k(i),-}(x)), \end{array}
    \right\}  \quad x\in [\alpha_i,\beta_i].
\end{align}
These relations reduce to
\begin{align}
    f_{i,+}(v) & = -f_{i,-}(v)), \qquad \textrm{for } v \in \ccal_{+,i}, \\
    f_{i,+}(v) & = f_{i,-}(v), \qquad \textrm{ for } v \in \ccal_{-,i}
\end{align}
and these jumps follow from \eqref{globalparametrix}, since we have chosen the
square root in $\sqrt{g(v)}$ with a cut along the union of arcs
$\bigcup_{i}\ccal_{+,i}$.

\subsection{Local parametrices around the endpoints: Airy parametrices}
\label{subsection:localparametrices}

In a small disk around the endpoints $\alpha_i$ and $\beta_i$ of the interval
$[\alpha_i,\beta_i]$ we construct  a local parametrix
$P^{(\textrm{Airy})}(z)$ involving Airy functions. Since the RH problem is
locally of size $2 \times 2$ and the equilibrium measures all vanish
as a square root, this can be done in the standard way \cite{Dei}. We
omit the details.

\subsection{Final transformation of the RH problem}
\label{subsection:errormatrix}

Using the global parametrix $P^{(\infty)}$ of
Section~\ref{subsection:globalparametrix} and the local parametrices
$P^{(\textrm{Airy})}$ of Section~\ref{subsection:localparametrices} we define
the final transformation $S \mapsto R$ of the RH problem by
\begin{equation}\label{defR}
R(z) = \left\{
    \begin{array}{ll}
    S(z)(P^{(\textrm{Airy})})^{-1}(z),& \quad \textrm{in the disks around } \alpha_i,\beta_i,\ i=1,\ldots,p+q-1, \\
    S(z)(P^{(\infty)})^{-1}(z),& \quad \textrm{elsewhere}.
\end{array}
\right.
\end{equation}

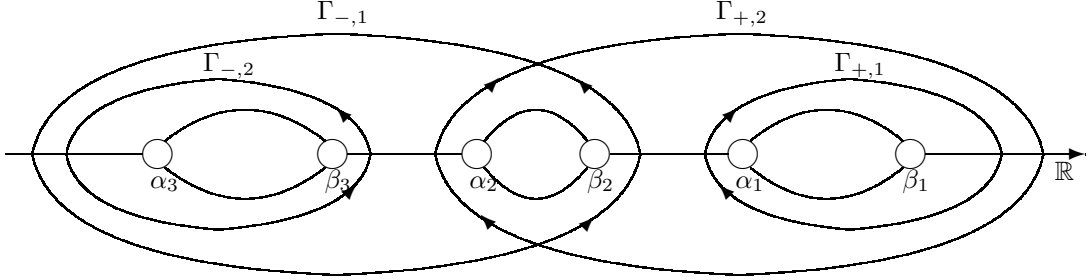
\begin{figure}[t]
\begin{center}
   \setlength{\unitlength}{1truemm}
   \begin{picture}(100,70)(-5,2)
       \qbezier(-14,40)(-12.5,31.5)(6,30)
       \qbezier(6,30)(24.5,31.5)(26,40)
       \qbezier(-14,40)(-12.5,48.5)(6,50)
       \qbezier(6,50)(24.5,48.5)(26,40)
       \put(4,51){$\Gamma_{-,2}$}
       \qbezier(34.5,40)(37.5,25.5)(74.5,24)
       \qbezier(74.5,24)(111.5,25.5)(114.5,40)
       \qbezier(34.5,40)(37.5,54.5)(74.5,56)
       \qbezier(74.5,56)(111.5,54.5)(114.5,40)
       \put(71.5,58){$\Gamma_{+,2}$}

       \put(23.2,35){\thicklines\vector(1,1){1}}
       \put(22.5,45.2){\thicklines\vector(-1,1){1}}
       \put(41.5,31){\thicklines\vector(-1,1){1}}
       \put(42,49.2){\thicklines\vector(1,1){1}}

       \qbezier(-18.5,40)(-15.5,25.5)(21.5,24)
       \qbezier(21.5,24)(58.5,25.5)(61.5,40)
       \qbezier(-18.5,40)(-15.5,54.5)(21.5,56)
       \qbezier(21.5,56)(58.5,54.5)(61.5,40)
       \put(19,58){$\Gamma_{-,1}$}
       \qbezier(70,40)(71.5,31.5)(89,30)
       \qbezier(89,30)(107.5,31.5)(109,40)
       \qbezier(70,40)(71.5,48.5)(89,50)
       \qbezier(89,50)(107.5,48.5)(109,40)
       \put(87,51){$\Gamma_{+,1}$}

       \put(72.8,35){\thicklines\vector(-1,1){1}}
       \put(73.3,45.2){\thicklines\vector(1,1){1}}
       \put(54.5,30.8){\thicklines\vector(1,1){1}}
       \put(54.2,49.1){\thicklines\vector(-1,1){1}}

       \put(116,36.6){$\mathbb R $}
       \put(120,40){\thicklines\vector(1,0){.0001}}
       \put(-3,35.6){$\alpha_{3}$}
       \put(20,35.6){$\beta_{3}$}
       \put(39,35.6){$\alpha_2$}
       \put(54.5,35.6){$\beta_2$}
       \put(74,35.6){$\alpha_1$}
       \put(96,35.6){$\beta_1$}

       \put(-2,40){\circle{4}}
       \put(21,40){\circle{4}}
       \put(40,40){\circle{4}}
       \put(55.5,40){\circle{4}}
       \put(75,40){\circle{4}}
       \put(97,40){\circle{4}}
       \qbezier(-1,38.2)(9.5,30)(20,38.2)
       \qbezier(-1,41.8)(9.5,50)(20,41.8)
       \qbezier(41,38.2)(47.75,30)(54.5,38.2)
       \qbezier(41,41.8)(47.75,50)(54.5,41.8)
       \qbezier(76,38.2)(86,30)(96,38.2)
       \qbezier(76,41.8)(86,50)(96,41.8)
       \put(-22,40){\line(1,0){18}}
       \put(23,40){\line(1,0){15}}
       \put(57.5,40){\line(1,0){15.5}}
       \put(99,40){\line(1,0){19}}
   \end{picture}
   \vspace{-20mm}
   \caption{The figure shows the contours in the RH problem for the final matrix $R(z)$ for the example
   of Figures \ref{fig:graph1} and \ref{fig:Riemannsurface4x4}.}
   \label{fig:finalRHP}
\end{center}
\end{figure}

From the construction of the parametrices it then follows that $R$ satisfies
the following RH problem.
\begin{rhp}\textrm{ }
\begin{itemize}
\item[\rm (1)] $R(z)$ is analytic in $\mathbb C \setminus \Sigma_R$
where $\Sigma_R$ is the contour shown in Figure~\ref{fig:finalRHP}.
\item[\rm (2)] $R$ has jumps $R_+ = R_- J_R$ on $\Sigma_R$,
where
\begin{align*}
    J_R(z) & = I_{p+q} + O(1/n), \qquad \text{on the boundaries of the disks}, \\
    J_R(z) & = I_{p+q} + O(e^{-cn(|z|+1)}), \qquad \text{on the other parts of $\Sigma_R$},
    \end{align*}
    for some constant $c >0$.
\item[\rm (3)] $R(z) = I_{p+q} + O(1/z)$ as $z \to \infty$.
\end{itemize}
\end{rhp}

As $n \to \infty$, the jump matrix $J_R$ tends to the identity matrix both in
$L^{\infty}(\Sigma_R)$ and in $L^{2}(\Sigma_R)$. Then as in
\cite{Dei,DKMVZ1,DKMVZ2} we may conclude that
\begin{equation} \label{asymptotics:R}
    R(z) = I_{p+q} + O\left(\frac{1}{n(|z|+1)}\right)
    \end{equation}
as $n \to \infty$, uniformly for $z$ in the complex plane. This completes the
RH steepest descent analysis.

\section{Proof of Theorem \ref{theorem:limitingdistribution}}
\label{section:proofmaintheorems}

Now we are ready to prove the main Theorem \ref{theorem:limitingdistribution}
by unfolding the transformations of the RH steepest descent analysis. Compare
with the proofs in the earlier papers \cite{BK2,BK3,DKV}.

For a finite $n$ we define the function $\rho^{(n)}$ as
\begin{align}
\label{def:rho}
    \rho^{(n)}(x) &= \frac{1}{\pi}\Im\xi^{(n)}_{k(i),+}(x), \qquad x\in [\alpha_i^{(n)},\beta_i^{(n)}],\quad i=1,\ldots,p+q-1,\\
    \rho^{(n)}(x) &= 0, \qquad \qquad\qquad\qquad x\in\mathbb R\setminus\bigcup_{i=1}^{p+q-1}[\alpha_i^{(n)},\beta_i^{(n)}].
\end{align}
Here we write $\xi^{(n)}_{k(i),+}(x)$ and $\alpha_i^{(n)}$, $\beta_i^{(n)}$ to emphasize
the $n$-dependence.

We recall that
\[ \Im\xi^{(n)}_{k(i),+} = -\Im\xi^{(n)}_{k(i),-} = -\Im\xi^{(n)}_{p+l(i),+} = \Im\xi^{(n)}_{p+l(i),-} \]
on the interval $[\alpha_i^{(n)},\beta_i^{(n)}]$, so one has in
fact several equivalent ways of expressing \eqref{def:rho}.

By \eqref{xifunctions:mixedangelesco:a}, \eqref{def:rho} and the
Stieltjes-Perron inversion formula one has that
\[ \rho^{(n)}(x) = \frac{d\mu_i^{(n)}(x)}{dx}, \qquad x\in [\alpha_i^{(n)},\beta_i^{(n)}]. \]

As $n \to \infty$, we have that $\alpha_i^{(n)} \to \alpha_i$, $\beta_i^{(n)} \to \beta_i$
and
\begin{equation} \label{rholimit}
    \lim_{n \to \infty} \rho^{(n)}(x) = \rho_i(x),
        \qquad x \in (\alpha_i, \beta_i),
    \end{equation}
where
\[ \rho_i(x) = \frac{d\mu_i(x)}{dx}, \qquad x \in [\alpha_i, \beta_i] \]
is the density of the $i$th component $\mu_i$ of the minimizer
$(\mu_1, \ldots, \mu_{p+q-1})$ of the vector
equilibrium problem with transition numbers $(t_{k,l})$.

Now we show that the $\rho_i$ give indeed the limiting distribution of the
non-intersecting Brownian motions. To this end we will use
\eqref{correlationkernellimit}. We start with the expression for the
correlation kernel \eqref{correlationkernel}, which we restate here for
convenience
\begin{equation*} \label{correlationkernel:Y}
    K_n(x,y) = \frac{1}{2\pi {\rm i}(x-y)}\begin{pmatrix} 0 & \cdots & 0 & w_{2,1}(y) &
    \cdots & w_{2,q}(y)\end{pmatrix} Y_{+}^{-1}(y)Y_{+}(x)\begin{pmatrix}
    w_{1,1}(x)\\ \vdots \\
    w_{1,p}(x)\\ 0 \\ \vdots \\ 0 \end{pmatrix}.
\end{equation*}
From the first transformation $Y\mapsto X$ in \eqref{defX}--\eqref{defL}
we get (we do not explicitly write the $n$-dependence in the $\lambda$-functions)
\begin{equation*} \label{correlationkernel:X}
    K_n(x,y) = \frac{1}{2\pi {\rm i}(x-y)}\begin{pmatrix} 0 & \cdots & 0 & e^{n\lambda_{p+1,+}(y)} &
    \cdots & e^{n\lambda_{p+q,+}(y)}\end{pmatrix} X_{+}^{-1}(y)X_{+}(x)\begin{pmatrix}
    e^{-n\lambda_{1,+}(x)}\\ \vdots \\
    e^{-n\lambda_{p,+}(x)}\\ 0 \\ \vdots \\ 0 \end{pmatrix}.
\end{equation*}
From the second transformation $X\mapsto T$ we obtain for
$x,y \in (\alpha_i^{(n)}, \beta_i^{(n)})$,
\begin{equation*} \label{correlationkernel:T}
    K_n(x,y) = \frac{1}{2\pi {\rm i}(x-y)}\left(e^{n\lambda_{p+l(i),+}(y)}\vece_{p+l(i)}^T\right) T_{+}^{-1}(y)T_{+}(x)
    \left(e^{-n\lambda_{k(i),+}(x)}\vece_{k(i)}\right).
\end{equation*}

From the third transformation $T\mapsto S$ in
\eqref{def:S} we get
\begin{multline} \label{correlationkernel:S}
    K_n(x,y) = \frac{1}{2\pi {\rm i}(x-y)}\left(-e^{n\lambda_{k(i),+}(y)}\vece_{k(i)}^T +
    e^{n\lambda_{p+l(i),+}(y)}\vece_{p+l(i)}^T
    \right) \\ \times S_{+}^{-1}(y)S_{+}(x)\left(e^{-n\lambda_{k(i),+}(x)}\vece_{k(i)}+
    e^{-n\lambda_{p+l(i),+}(x)}\vece_{p+l(i)}\right),
\end{multline}
for $x,y\in (\alpha_i^{(n)},\beta_i^{(n)})$.
Defining the function $h_n$ by
\begin{equation}\label{def:h} h_n(x) := -\Re(\lambda_{k(i),+}(x)) =
-\Re(\lambda_{p+l(i),+}(x)),\quad x\in [\alpha_i^{(n)},\beta_i^{(n)}],
\end{equation}
we see that \eqref{correlationkernel:S} can be rewritten as
\begin{multline} \label{correlationkernel:Sbis}
    K_n(x,y) = \frac{e^{n(h_n(x)-h_n(y))}}{2\pi {\rm i}(x-y)}\left(-e^{n{\rm i}\Im(\lambda_{k(i),+}(y))}\vece_{k(i)}^T +
    e^{-n{\rm i}\Im(\lambda_{k(i),+}(y))}\vece_{p+l(i)}^T
    \right) \\ \times S_{+}^{-1}(y)S_{+}(x)\left(e^{-n{\rm i}\Im(\lambda_{k(i),+}(x))}\vece_{k(i)}+
    e^{n{\rm i}\Im(\lambda_{k(i),+}(x))}\vece_{p+l(i)}\right),
\end{multline}
for $x,y\in (\alpha_i^{(n)},\beta_i^{(n)})$.

Now from \eqref{asymptotics:R} it follows by standard arguments (e.g.\
\cite[Section 9]{BK2}) that
\begin{equation*}
S^{-1}(y)S(x) = I+O(x-y),\quad \textrm{as }y\to x
\end{equation*}
uniformly in $n$. Hence \eqref{correlationkernel:Sbis} takes the form
\begin{equation} \label{correlationkernel:R}
    K_n(x,y) = e^{n(h_n(x)-h_n(y))}
    \left( \frac{\sin(n\Im(\lambda_{k(i),+}(x)-\lambda_{k(i),+}(y)))}{\pi(x-y)} +O(1)
    \right),
\end{equation}
for $x,y\in (\alpha_i,\beta_i)$, where the $O(1)$ term holds uniformly in $n$. Then by letting
$y \to x$ and using l'H\^opital's rule we find
\begin{equation*}
    K_n(x,x) = \frac{n}{\pi}(\Im(\xi^{(n)}_{k(i),+}(x)))+O(1),
\end{equation*}
for $x\in (\alpha_i,\beta_i)$, or equivalently
\begin{equation*}
    K_n(x,x) = n\rho^{n}(x)+O(1),
\end{equation*}
by virtue of \eqref{def:rho}. It follows from \eqref{rholimit} that
\begin{equation} \label{correlationkernel:Rbis}
    \lim_{n\to\infty} \frac{1}{n} K_n(x,x) = \rho_i(x),\qquad x \in (\alpha_i,\beta_i).
\end{equation}
In a similar way one can prove that  \begin{equation}
\label{correlationkernel:Rtris}
    \lim_{n\to\infty}\frac{1}{n} K_n(x,x) = 0,\quad x\in\mathbb R\setminus\bigcup_i [\alpha_i,\beta_i].
\end{equation}
This completes the proof of Theorem \ref{theorem:limitingdistribution}.


\begin{thebibliography}{99}

\bibitem{AFvM}
    M. Adler, P.L. Ferrari, and P. van Moerbeke,
    Airy processes with wanderers and new universality classes,
    preprint math-pr/0811.1863.
\bibitem{Apt}
    A.I. Aptekarev,
    Multiple orthogonal polynomials,
    J. Comput. Appl. Math. 99 (1998), 423–-447.
\bibitem{ABK}
    A.I. Aptekarev, P.M. Bleher, and A.B.J. Kuijlaars,
    Large $n$ limit of Gaussian random matrices with external
    source, part II, Comm. Math. Phys. 259 (2005), 367--389.
\bibitem{BK}
    P.M. Bleher and A.B.J. Kuijlaars,
    Random matrices with external source and multiple orthogonal polynomials,
    Int. Math. Res. Not. 2004, no.~3 (2004), 109–-129.
\bibitem{BK2}
    P.M. Bleher and A.B.J. Kuijlaars,
    Large $n$ limit of Gaussian random matrices with external
    source, part I, Comm. Math. Phys. 252 (2004), 43--76.
\bibitem{BK3}
    P.M. Bleher and A.B.J. Kuijlaars,
    Large $n$ limit of Gaussian random matrices with external
    source, part III: double scaling limit,
    Comm. Math. Phys. 270 (2007), 481--517.
\bibitem{Bor}
    A. Borodin,
    Borthogonal ensembles,
    Nucl. Phys. B536 (1999), 704--732.
\bibitem{Claeys1}
    T. Claeys and A.B.J. Kuijlaars,
    Universality of the double scaling limit in random matrix
    models,
    Comm. Pure Appl. Math. 59 (2006), 1573--1603.
\bibitem{DK2}
    E. Daems and A.B.J. Kuijlaars,
    Multiple orthogonal polynomials of mixed type and non-intersecting Brownian
    motions, J. Approx. Theory 146 (2007), 91--114.
\bibitem{DKV}
    E. Daems, A.B.J. Kuijlaars, and W. Veys,
    Asymptotics of non-intersecting Brownian motions and a $4\times 4$ Riemann-Hilbert problem,
    J. Approx. Theory. 153 (2008), 225--256.
\bibitem{Dei}
    P. Deift, Orthogonal Polynomials and Random Matrices: a Riemann-Hilbert
    approach. Courant Lecture Notes in Mathematics Vol. 3, Amer. Math. Soc.,
    Providence R.I. 1999.
\bibitem{DKM}
    P. Deift, T. Kriecherbauer, and K.T-R McLaughlin,
    New results on the equilibrium measure for logarithmic potentials in the presence of an external
    field,
    J. Approx. Theory 95 (1998), 388--475.
\bibitem{DKMVZ1}
    P. Deift, T. Kriecherbauer, K.T-R McLaughlin, S. Venakides, and  X. Zhou,
    Uniform asymptotics for polynomials orthogonal with respect to
    varying exponential weights and applications to universality
    questions in random matrix theory,
    Comm. Pure Appl. Math. 52 (1999), 1335--1425.
\bibitem{DKMVZ2}
    P. Deift,  T. Kriecherbauer, K.T-R  McLaughlin, S. Venakides, and  X. Zhou,
    Strong asymptotics of orthogonal polynomials with respect to exponential weights,
    Comm. Pure Appl. Math. 52 (1999), 1491--1552.
\bibitem{Del}
    S. Delvaux,
    Average characteristic polynomials for multiple orthogonal polynomial ensembles,
    preprint math.CA/0907.0156.
\bibitem{DelKui}
    S. Delvaux and A.B.J. Kuijlaars,
    A phase transition for non-intersecting Brownian motions, and the Painlev\'e
    II equation, preprint math.CV/08091000. To appear in Int. Math. Res. Not.
\bibitem{Duits}
    M. Duits and A.B.J. Kuijlaars,
    Universality in the two matrix model: A Riemann-Hilbert steepest descent
    analysis, Comm. Pure Appl. Math. 62 (2009), 1076--1153.
\bibitem{Dyson}
    F.J. Dyson,
    A Brownian-motion model for the eigenvalues of a random matrix,
    J. Math. Phys. 3 (1962), 1191--1198.
\bibitem{FIK}
    A.S. Fokas, A.R. Its, and A.V. Kitaev,
    The isomonodromy approach to matrix models in 2D quantum gravity,
    Commun. Math. Phys. 147 (1992), 395--430.
\bibitem{GVL}
    G.H. Golub and C.F. Van Loan,
    Matrix Computations, The Johns
    Hopkins University Press, third edition, 1996.
\bibitem{GRS}
    A.A. Gonchar, E.A. Rakhmanov, and V.N. Sorokin,
    Hermite-Pad\'e approximants for systems of Markov-type functions,
    Sbornik Mathematics 188 (1997), 671--696.
\bibitem{Joh1}
    K. Johansson,
    Random matrices and determinantal processes, in: Mathematical Statistical
    Physics: Lecture Notes of the Les Houches Summer School 2005 (Bovier et al., eds.),
    Elsevier, 2006, pp. 1-–55.
\bibitem{KMcG}
    S. Karlin and J. McGregor, Coincidence probabilities,
    Pacific J. Math., 9 (1959), 1141--1164.
\bibitem{KKPS}
    A.M. Khorunzhy, B.A. Khoruzhenko, L.A. Pastur, and M.V. Shcherbina,
    The large $n$-limit in statistical mechanics and the spectral theory of disordered
    systems. In: Phase transitions and critical phenomena, 15 (C. Domb and J.~L. Lebowitz eds),
    Academic Press, London, 1992, pp.~73--239.
\bibitem{Mcl}
    K. T.-R. McLaughlin,
    Asymptotic analysis of random matrices with external source and a family of
    algebraic curves,  Nonlinearity 20 (2007), 1547--1571.
\bibitem{Mir}
    R. Miranda,
    Algebraic Curves and Riemann Surfaces,
    American Mathematical Society, Providence, RI, 1995.
\bibitem{NS}
    E. M. Nikishin and V. N. Sorokin,
    Rational Approximations and Orthogonality,
    Amer. Math. Soc., Providence, RI, 1991.
\bibitem{Ora}
    N. Orantin,
    Gaussian matrix model in an external field and non-intersecting Brownian motions,
    preprint math-ph/0803.0705.
\bibitem{Pastur}
    L. Pastur,
    The spectrum of random matrices (Russian),
    Teoret. Mat. Fiz. 10 (1972), 102--112.
\bibitem{SaffTotik}
    E.B. Saff and V. Totik,
    Logarithmic Potentials with External Field,
    Springer-Verlag, Berlin, 1997.
\bibitem{Totik}
    V. Totik,
    Weighted approximation with varying weight,
    Lecture Notes in Mathematics 1569,
    Springer-Verlag, Berlin, 1994.
\bibitem{VAGK} W. Van Assche, J.S. Geronimo, and A.B.J. Kuijlaars,
    Riemann-Hilbert problems for multiple orthogonal polynomials,
     Special Functions 2000: Current Perspectives and Future Directions
    (J. Bustoz et al., eds.), Kluwer, Dordrecht, 2001, pp. 23--59.
\end{thebibliography}
\end{document}